\documentstyle[twocolumn,eqsecnum,aps]{revtex}

\begin{document}

\draft

\title{Convergence Acceleration via Combined Nonlinear-Condensation
Transformations}

\author{Ulrich D.\ Jentschura\thanks{Internet: ulj@nist.gov}}
\address{National Institute of Standards and Technology, Gaithersburg,
MD 20899-0001, USA \\
and \\
Institut f\"{u}r Theoretische Physik, TU Dresden, D-01062 Dresden,
Germany}

\author{Peter J.\ Mohr\thanks{Internet: mohr@nist.gov}}
\address{National Institute of Standards and Technology, Gaithersburg,
MD 20899-0001, USA}

\author{Gerhard Soff\thanks{Internet: soff@physik.tu-dresden.de}}
\address{Institut f\"{u}r Theoretische Physik, TU Dresden, D-01062
Dresden, Germany}

\author{Ernst Joachim
Weniger\thanks{Internet: joachim.weniger@chemie.uni-regensburg.de}}
\address{Center for Theoretical Studies of Physical Systems, Clark
Atlanta University, James P.\ Brawley Drive at Fair Street, S.W.,
Atlanta, GA 30314, USA \\
and \\
Institut f\"{u}r Physikalische und Theoretische Chemie, Universit\"{a}t
Regensburg, D-93040 Regensburg, Germany}

\date{Submitted to Computer Physics Communications}
\maketitle

\begin{abstract}
A method of numerically evaluating slowly convergent monotone series is
described. First, we apply a condensation transformation due to Van
Wijngaarden to the original series. This transforms the original
monotone series into an alternating series. In the second step, the
convergence of the transformed series is accelerated with the help of
suitable nonlinear sequence transformations that are known to be
particularly powerful for alternating series. Some theoretical aspects
of our approach are discussed. The efficiency, numerical stability, and
wide applicability of the combined nonlinear-condensation transformation
is illustrated by a number of examples. We discuss the evaluation of
special functions close to or on the boundary of the circle of
convergence, even in the vicinity of singularities. We also consider a
series of products of spherical Bessel functions, which serves as a
model for partial wave expansions occurring in quantum electrodynamic
bound state calculations.
\end{abstract}


\vspace{0.05in}
\vspace{0.05in}


\section{Introduction}

Divergent and slowly convergent series occur abundantly in the
mathematical and physical sciences. Accordingly, there is an extensive
literature on numerical techniques which convert a divergent or slowly
convergent series into a new series with hopefully better numerical
properties. An overview of the existing sequence transformations as well
as many references can be found in books by Wimp \cite{Wi81} and
Brezinski and Redivo Zaglia \cite{BreRZa91}. The historical development
of these techniques up to 1945 is described in a monograph by Brezinski
\cite{Br91}, and the more recent developments are discussed in an
article by Brezinski \cite{Bre96}.

A very important class of sequences $\{ s_n\}_{n=0}^{\infty}$ is
characterized by the asymptotic condition
\begin{equation}
\lim_{n \to \infty} \frac {s_{n+1} - s} {s_n - s} \; = \; \rho \, ,
\label{LimRatSeq}
\end{equation}
which closely resembles the well known ratio test for infinite series.
Here, $s = s_{\infty}$ is the limit of this sequence as $n \to \infty$.
A convergent sequence satisfying condition (\ref{LimRatSeq}) with $\vert
\rho \vert < 1$ is called {\em linearly} convergent, and it is called
{\em logarithmically} convergent if $\rho = 1$.

If the elements of the sequence $\{ s_n \}_{n=0}^{\infty}$ in Eq.\
(\ref{LimRatSeq}) are the partial sums $s_n = \sum_{k=0}^{n} a_k$ of an
infinite series, and if $\rho$ satisfies either $\rho = 1$ or $0 < \rho
< 1$ (these are the only cases which will be considered in this
article), then there exists an integer $N \ge 0$ such that all terms
$a_k$ with $k \ge N$ have the same sign. Hence, series of this kind can
be split up into a finite sum containing the leading terms with the
irregular signs, and a monotone series whose terms all have the same
sign. The subject of this article is the efficient and reliable
evaluation of series that exhibit these properties.

The partial sums of a nonterminating Gaussian hypergeometric series
\begin{equation}
{}_2 F_1 (a, b; c; z) \; = \;
\sum_{m=0}^{\infty} \, \frac {(a)_m (b)_m} {(c)_m} \, \frac{z^m}{m!} \, ,
\label{2F1_Ser}
\end{equation}
where $(a)_m = \Gamma (a+m)/\Gamma (a) = a (a+1)\ldots (a+m-1)$ is a
Pochhammer symbol, or its generalization
\begin{eqnarray}
\lefteqn{{}_{p+1} F_p \bigl(\alpha_1, \ldots, \alpha_{p+1};
\beta_1, \ldots, \beta_p; z \bigr)} \nonumber \\
& \qquad = &
\sum_{m=0}^{\infty} \, \frac
{(\alpha_1)_m\ldots(\alpha_{p+1})_m}
{(\beta_1)_m\ldots(\beta_p)_m} \, \frac{z^m}{m!} \, ,
\label{GenHygSer}
\end{eqnarray}
which converge for $\vert z \vert < 1$ and diverge for $\vert z \vert >
1$, are typical examples of linearly convergent sequences with $\rho =
z$.

The partial sums of the Dirichlet series
\begin{equation}
\zeta (z) \; = \; \sum_{m=0}^{\infty} \, (m+1)^{- z}
\label{ZetaFun}
\end{equation}
for the Riemann zeta function converge logarithmically if ${\rm Re} (z)
> 1$. This follows from the following asymptotic estimate ($n
\to \infty$) of the truncation error (p.\ 21 of Ref.\ \cite{Wi81}):
\begin{eqnarray}
\lefteqn{\zeta (z) \, - \, \sum_{m=0}^{n} \, (m+1)^{- z}} \nonumber \\
& \qquad = & \frac {(n+1)^{1-z}} {z-1} \, - \, \frac {1} {2 (n+1)^z}
\, + \, O \bigl(n^{- z - 1}\bigr) \, .
\label{ZetaErrEst}
\end{eqnarray}
The Dirichlet series (\ref{ZetaFun}) is notorious for its extremely slow
convergence if ${\rm Re} (z)$ is only slightly larger than one. In this
case, the series can only be used for the computation of $\zeta (z)$ if
it is combined with suitable convergence acceleration methods like the
Euler-Maclaurin summation (see for instance Chapter 8 of Ref.\ 
\cite{Olv74}, p.\ 379 of Ref.\ \cite{BeOr78}, or Chapter 6 of Ref.\ 
\cite{Edw74}).

The acceleration of linear convergence is comparatively simple, both
theoretically and practically, as long as $\rho$ in Eq.\
(\ref{LimRatSeq}) is not too close to one. With the help of
Germain-Bonne's formal theory of convergence acceleration \cite{Ge-Bo73}
and its extension (Section 12 of Ref.\ \cite{We89}), it can be decided
whether a sequence transformation is capable of accelerating linear
convergence or not. Moreover, many sequence transformations are known
which are capable of accelerating linear convergence {\em effectively}.
Examples are the $\Delta^2$ process, which is usually attributed to
Aitken \cite{Ai26} although it is in fact much older (p.\ 90 of Ref.\
\cite{Br91}), Wynn's epsilon algorithm \cite{Wy56a}, which produces
Pad\'{e} approximants if the input data are the partial sums of a formal
power series, or Levin's sequence transformation \cite{Le73} and
generalizations (Sections 7 - 9 of Ref.\ \cite{We89}), which require as
input data not only the elements of the sequence to be transformed but
also explicit estimates for the corresponding truncation errors.

The acceleration of logarithmic convergence is much more difficult than
the acceleration of linear convergence. Delahaye and Germain-Bonne
\cite{DeGe-Bo82} showed that no sequence transformation can exist which
is able to accelerate the convergence of {\em all\/} logarithmically
convergent sequences. Consequently, in the case of logarithmic
convergence the success of a convergence acceleration process cannot be
guaranteed unless additional information is available. Also, an analogue
of Germain-Bonne's formal theory of linear convergence acceleration
\cite{Ge-Bo73} cannot exist.

In spite of these complications, many sequence transformations are known
which work reasonably well at least for suitably restricted subsets of
the class of logarithmically convergent sequences. Examples are
Richardson extrapolation \cite{Ri27}, Wynn's rho algorithm \cite{Wy56b}
and its iteration (Section 6 of Ref.\ \cite{We89}) as well as Osada's
modification of the rho algorithm \cite{Os90}, Brezinski's theta algorithm
\cite{Bre71} and its iteration (Section 10 of Ref.\ \cite{We89}),
Levin's $u$ and $v$ transformations \cite{Le73} and related
transformations (Sections 7 - 9 of Ref.\ \cite{We89}), and the
modification of the $\Delta^2$ process by Bj{\o}rstad, Dahlquist, and
Grosse \cite{BjDaGr81}. However, there is a considerable amount of
theoretical and empirical evidence that sequence transformations are in
general less effective in the case of logarithmic convergence than in
the case of linear convergence.

Numerical stability is a very important issue. A sequence transformation
can only accelerate convergence if it succeeds in extracting some
additional information about the index-dependence of the truncation
errors from a finite set $s_n$, $s_{n+1}$, $\ldots$, $s_{n+k}$ of input
data. Normally, this is done by forming arithmetic expressions involving
higher weighted differences. However, forming higher weighted
differences is a potentially unstable process which can easily lead
to a serious loss of significant digits or even to completely
nonsensical results.

If the input data are the partial sums of a {\em strictly alternating\/}
series, the formation of higher weighted differences is normally a
remarkably stable process. Hence, a serious loss of significant digits
is not to be expected if the partial sums of a strictly alternating
series are used as input data in a convergence acceleration or summation
process. If, however, the input data are the partial sums of a {\em
monotone\/} series, numerical instabilities due to cancellation are
quite likely, in particular if convergence is very slow. Thus, if the
sequence to be transformed either converges linearly with a value of
$\rho$ in Eq.\ (\ref{LimRatSeq}) that is only slightly smaller than one,
or if it converges logarithmically ($\rho = 1$), numerical instabilities
are a serious problem and at least some loss of significant digits is to
be expected.

Generally, the sequence transformations mentioned above are not able to
determine the limit of a logarithmically convergent series, whose terms
ultimately all have the same sign, with an accuracy close to machine
accuracy. This restricts the practical usefulness of sequence
transformations severely, e.g., in FORTRAN calculations with a fixed
precision.

In this paper, we show that these stability problems can often be
overcome by transforming slowly convergent monotone series not by
straightforward application of a {\em single\/} sequence transformation
but by a combination of two {\em different\/} transformations.

In the {\em first\/} step, a monotone series is transformed into a
strictly alternating series with the help of a {\em condensation
transformation}. This transformation was first mentioned on p.\ 126 of
Ref.\ \cite{NaPhyLa61} and only later published by Van Wijngaarden
\cite{vWij65}. Later, the Van Wijngaarden transformation was studied by
Daniel \cite{Dan69}, and recently it was rederived by Pelzl and King
\cite{PeKi98}, who used it for the high-precision calculation of atomic
three-electron interaction integrals of explicitly correlated wave
functions.

In the {\em second\/} step, the convergence of the resulting alternating
series is accelerated by suitable {\em nonlinear sequence
transformations} \cite{We89,Le73} which are known to be very powerful in
the case of alternating series. Since the transformation of alternating
series is a remarkably stable process, the limits of even extremely
slowly convergent monotone series can be determined with an accuracy
close to machine accuracy. Conceptually, but not technically our
approach resembles that of Brezinski, Delahaye, and Germain-Bonne
\cite{BreDelGe-Bo83} who proposed to extract a linearly convergent
subsequence from a logarithmically convergent input sequence by a
selection process.

In this article, we will call our approach, which consists of the Van
Wijngaarden {\em condensation\/} transformation and the subsequent {\em
nonlinear\/} sequence transformation, the {\em combined
nonlinear-condensation transformation} (CNCT).

The CNCT is not restricted to logarithmically convergent series. It can
also be used in the case of a linearly convergent monotone series with a
value of $\rho$ in Eq.\ (\ref{LimRatSeq}) close to one. Typically, this
corresponds to a power series whose coefficients ultimately all have the
same sign and whose argument is positive and close to the boundary of
the circle of convergence. We will present some examples which show that
the CNCT works even if the argument of the power series is very close to
a singularity.

Our approach requires the evaluation of terms of the original monotone
series with high indices. Consequently, our two-step approach is
computationally more demanding than the application of a single sequence
transformation, and it cannot be applied if only a few terms of a slowly
convergent series are available. In spite of these restrictions, we
believe that the CNCT is very useful at least for certain problems since
it is able to produce highly accurate results that can only be
accomplished otherwise with a considerably greater numerical effort.

Special functions are defined and, in many cases, also evaluated via
series expansions. The evaluation of special functions is an old problem
of numerical mathematics with a very extensive literature (compare for
example the books by Luke \cite{Luke69,Luke77}, Van Der Laan and Temme
\cite{vdLaTe80}, and Zhang and Jin \cite{ZhaJin96}). Nevertheless,
there is still a considerable amount of research going on and many new
algorithms for the computation of special functions have been developed
recently (compare for example the papers by Lozier and Olver
\cite{LozOlv48} and Lozier \cite{Loz96} and the long lists of references
therein). We believe that sequence transformations in general and the
CNCT in particular are useful tools for the evaluation of special
functions. Of course, this is also true for problems in theoretical and
computational physics.

This paper is organized as follows: In Section \ref{sec:VWT}, we
describe the Van Wijngaarden transformation. In Section
\ref{sec:NonSeqTrans}, we discuss the nonlinear sequence transformations
which we use for the acceleration of the resulting alternating series.
In Section \ref{sec:ZetaFun}, we apply the CNCT to the Riemann zeta
function. In Section \ref{sec:LerTran}, we consider the evaluation of
the Lerch transcendent and related functions with arguments on or close
to the boundary of the circle of convergence. In Section
\ref{sec:GenHygFun}, we examine the evaluation of the generalized
hypergeometric series ${}_{p+1} F_p$ ($p \ge 2$) with unit argument or
with an argument $z$ which is only slightly smaller than one. In Section
\ref{sec:ProdBessHank}, we discuss the evaluation of an infinite series
involving Bessel and Hankel functions. In Appendix \ref{App_A}, we
discuss the efficiency of sequence transformations in the case of
monotone series or strictly alternating series in more detail. Finally,
in Appendix \ref{App_B} we discuss exactness properties of the sequence
transformations which we use in the second step for the acceleration of
the convergence of the resulting alternating series.

The example of Section \ref{sec:ProdBessHank} serves as a model problem
for slowly convergent series, which occur in quantum electrodynamic
bound state calculations and which were treated successfully with the
methods discussed in this paper \cite{JeSoMo98}. Therefore, we expect
the CNCT to be a general computational tool for the evaluation of slowly
convergent sums over intermediate angular momenta which arise from the
decomposition of relativistic propagators in QED bound state
calculations into partial waves.

All calculations were done in Mathematica\footnote{Certain commercial
equipment, instruments, or materials are identified in this paper to
foster understanding. Such identification does not imply recommendation
or endorsement by the National Institute of Standards and Technology,
nor does it imply that the materials or equipment identified are
necessarily the best available for the purpose.} with a relative
accuracy of 16 decimal digits \cite{mathematica}. In this way, we
simulate the usual DOUBLE PRECISION accuracy in FORTRAN.

\section{The Van Wijngaarden transformation}
\label{sec:VWT}

Let us assume that the partial sums
\begin{equation}
\sigma_n \; = \; \sum_{k=0}^{n} \, a (k)
\label{ParSumInput}
\end{equation}
of an infinite series converge either linearly or logarithmically to
some limit $\sigma = \sigma_{\infty}$ as $n \to \infty$. We also assume
that all terms $a (k)$ have the same sign, i.e., the series
$\sum_{k=0}^{\infty} a (k)$ is a monotone series.

Following Van Wijngaarden \cite{vWij65}, we transform the original
series into an alternating series $\sum_{j=0}^{\infty} (-1)^j {\bf A}_j$
according to
\begin{eqnarray}
\sum_{k=0}^{\infty} \, a (k) & = &
\sum_{j=0}^{\infty} \, (-1)^j \, {\bf A}_j \, ,
\label{VWSerTran}
\\
{\bf A}_j & = & \sum_{k=0}^{\infty} \, {\bf b}_{k}^{(j)} \, ,
\label{A2B}
\\
{\bf b}_{k}^{(j)} & = & 2^k \, a(2^k\,(j+1)-1) \, .
\label{B2a}
\end{eqnarray}
Obviously, the terms ${\bf A}_j$ defined in Eq.\ (\ref{A2B}) all have
the same sign if the terms $a (k)$ of the original series all have the
same sign. In the sequel, the quantities ${\bf A}_j$ will be referred to
as the {\em condensed series\/}, and the series $\sum_{j=0}^{\infty}
(-1)^j {\bf A}_j$ will be referred to as the {\em transformed
alternating series}, or alternatively as the {\em Van Wijngaarden
transformed series}.

We call the Van Wijngaarden transformation a condensation transformation
because its close connection to Cauchy's condensation theorem (p.\ 28 of
Ref.\ \cite{Bro91} or p.\ 121 of Ref.\ \cite{Kn64}). Given a monotone
series $\sum_{k=0}^{\infty} a (k)$ with terms that satisfy $\vert a(k+1)
\vert < \vert a(k) \vert$, Cauchy's condensation theorem states that
$\sum_{k=0}^{\infty} a (k)$ converges if and only if the first condensed
series ${\bf A}_0$ defined in Eq.\ (\ref{A2B}) converges.

In Eq.\ (\ref{A2B}), the indices of the terms of the original series are
chosen in such a way that sampling at very high indices takes place
(according to Eq.\ (\ref{B2a}), the indices of the terms of the original
series grow exponentially). In this way, we obtain information about the
behaviour of the terms of the original series at high indices.

Moreover, if the terms of the original series behave asymptotically ($n
\to \infty$) either like $a (n) \sim n^{-1 - \epsilon}$ with $\epsilon >
0$ or like $a (n) \sim n^{\beta} r^n$ with $0 < r < 1$ and $\beta$ real,
then the terms of the original series become negligibly small after a
few evaluations. Specifically, the series (\ref{A2B}) for the terms
${\bf A}_j$ converges linearly in these cases.

When summing over $k$ in Eq.\ (\ref{A2B}), we found that it is normally
sufficient to terminate this sum when the last term is smaller than the
desired accuracy. Typically, 20 to 30 terms are needed for
a relative accuracy of $10^{-14}$ in the final result for
the condensed sum.

It should be possible to accelerate the convergence of the series
(\ref{A2B}) for ${\bf A}_j$. Since, however, this is a monotone series,
it is not clear whether and how many digits would be lost in the
convergence acceleration process. Consequently, we prefer to perform 
a safe and straightforward evaluation of the condensed sums ${\bf A}_j$
and add up the terms of the series until convergence is reached,
although this is most likely not the most efficient approach.

The transformation from a monotone series to a strictly alternating
series according to Eqs.\ (\ref{VWSerTran}) - (\ref{B2a}) is essentially
a reordering of the terms $a(k)$ of the original series. This is seen as
follows. We first define the partial sums
\begin{equation}
{\bf S}_n \; = \; \sum_{j=0}^{n} \, (-1)^j \, {\bf A}_j
\label{PsumS}
\end{equation}
of the Van Wijngaarden transformed original series. It can be shown
easily that ${\bf S}_n$ with $n \ge 0$ reproduces the partial sum
$\sigma_n$, Eq. (\ref{ParSumInput}), which contains the first $n+1$
terms of the original series. To illustrate this procedure, we present
Table \ref{Tab_2_1}. Formal proofs of the correctness of this
rearrangement can be found in Ref.\ \cite{Dan69} or in the Appendix of
Ref.\ \cite{PeKi98}.

\medskip
\begin{center}
{\bf \Large Insert Table \ref{Tab_2_1} here}
\end{center}
\medskip

Daniel \cite{Dan69} was able to formulate some mild conditions which
guarantee that the limits $\sigma$ and ${\bf S}$ of the partial sums
(\ref{ParSumInput}) and (\ref{PsumS}), respectively, simultaneously
exist and are equal:
\begin{equation}
\lim_{n \to \infty} \, \sigma_n \; = \; \sigma \; = \;
\lim_{n \to \infty} \, {\bf S}_n \; = \; {\bf S} \, .
\end{equation}
For example, in the Corollary on p.\ 92 of Ref.\ \cite{Dan69} it was
shown that if a strictly decreasing sequence $\{ M_k \}_{k=0}^{\infty}$
of positive bounds exists which satisfy $\vert a(k) \vert \le M_k$ for
all $k \ge 0$ {\em and} if $\sum_{k=0}^{\infty} M_k < \infty$ holds,
then the original monotone series and the Van Wijngaarden transformed
series on the right-hand side of Eq.\ (\ref{VWSerTran}) both converge to
the same limit (i.e., $\sigma = {\bf S}$ holds). This useful criterion
is fulfilled by all series considered in this paper.

\section{Nonlinear Sequence Transformations}
\label{sec:NonSeqTrans}

The series (\ref{A2B}) for the terms of the Van Wijngaarden transformed
series can be rewritten as follows:
\begin{equation}
{\bf A}_j \; = \; a (j) \, + \, 2 a (2j+1) \, + \, 4 a (4j+3) \, + \,
\ldots \, .
\end{equation}
Since the terms $a (k)$ of the original series have by assumption the
same sign, we immediately observe
\begin{equation}
\vert {\bf A}_j \vert \ge \vert a_j \vert \, .
\label{aBoundA}
\end{equation}
Consequently, the Van Wijngaarden transformation, Eqs.\
(\ref{VWSerTran}) - (\ref{B2a}), does not lead to an alternating series
whose terms decay more rapidly in magnitude than the terms of the
original monotone series. Thus, an acceleration of convergence can only
be achieved if the partial sums (\ref{PsumS}) of the Van Wijngaarden
transformed series are used as input data in a convergence acceleration
process.

Since the Van Wijngaarden transformed series on the right-hand side of
Eq.\ (\ref{VWSerTran}) is alternating if the terms of the original
series all have the same sign, it is recommended to choose a suitable
convergence acceleration method which is particularly powerful in the
case of alternating series. A judicious choice of the convergence
accelerator is of utmost importance since it ultimately decides whether
our approach is numerically useful or not.

Daniel \cite{Dan69} used the Euler transformation (Eq.\ (3.6.27) on p.\
16 of Ref.\ \cite{AbrSte72})):
\begin{equation}
\sum_{k=0}^{\infty} \, (-1)^k u_k \; = \; \sum_{k=0}^{\infty} \,
\frac{(-1)^k}{2^{k+1}} \, \Delta^k u_0 \, .
\label{EuTrans}
\end{equation}
Here, $\Delta$ is the (forward) difference operator defined by
$\Delta f (n) = f (n+1) - f (n)$, and
\begin{equation}
\Delta^k u_0 \; = \; (-1)^k \,
\sum_{m=0}^{k} \, (-1)^m \, {{k} \choose {m}} \, u_{m} \, .
\end{equation}
The Euler transformation, which was published in its original version in
1755 on p.\ 281 of Ref. \cite{Euler1755}, is a series transformation
which is specially designed for alternating series. It is treated in
many books on numerical mathematics. Nevertheless, the Euler
transformation is not a particularly efficient accelerator for the Van
Wijngaarden transformed series considered in this article. In Table
\ref{Tab_4_1}, we show some explicit results obtained by the Euler
transformation. We also applied the Euler transformation to all other
Van Wijngaarden transformed series presented in this paper, and we
consistently observed that it is clearly less powerful than the
transformations we used.

Much better results can be expected from the more modern nonlinear
sequence transformations which transform a sequence $\{ s_n
\}_{n=0}^{\infty}$, whose elements may be the partial sums of an
infinite series, into a new sequence $\{ s^{\prime}_n \}_{n=0}^{\infty}$
with better numerical properties \cite{Wi81,BreRZa91,We89}.

The basic assumption of a sequence transformation is that the elements
of the sequence $\{ s_n \}_{n=0}^{\infty}$ to be transformed can be
partitioned into a limit $s$ and a remainder $r_n$ according to
\begin{equation}
s_n \; = \; s \, + \, r_n 
\end{equation}
for all $n \ge 0$. Only in the case of some more or less trivial model
cases will a sequence transformation be able to determine the limit $s$
of the sequence $\{ s_n \}_{n=0}^{\infty}$ exactly after a {\em
finite\/} number of steps. Hence, the elements of the transformed
sequence $\{ s^{\prime}_n \}_{n=0}^{\infty}$ can be partitioned into the
same limit $s$ and a transformed remainder $r_{n}^{\prime}$ according to
\begin{equation}
s_{n}^{\prime} \; = \; s \, + \, r_{n}^{\prime} 
\end{equation}
for all $n \ge 0$. In general, the transformed remainders
$r_{n}^{\prime}$ are nonzero for all finite values of $n$. However,
convergence is obviously increased if the transformed remainders $\{
r^{\prime}_n \}_{n=0}^{\infty}$ vanish more rapidly than the original
remainders $\{ r_n \}_{n=0}^{\infty}$:
\begin{equation}
\lim_{n \to \infty} \, \frac {r_{n+1}}{r_n} \; = \;
\lim_{n \to \infty} \, \frac {s_{n+1}-s}{s_n-s} \; = \; 0 \, .
\end{equation}

The best known nonlinear sequence transformation is probably Wynn's
epsilon algorithm \cite{Wy56a} which produces Pad\'{e} approximants if
the input data are the partial sums of a formal power series. The
epsilon algorithm is also known to work well in the case of alternating
series. Consequently, it is an obvious idea to use the epsilon algorithm
for the acceleration of the alternating series which we obtain via the
Van Wijngaarden transformation.

However, much better results can be obtained by applying sequence
transformations which use explicit remainder estimates as input data in
addition to the partial sums of the series to be transformed.
Consequently, in this article we use exclusively the following two
sequence transformations by Levin \cite{Le73} and Weniger \cite{We89}, 
respectively:
\begin{eqnarray}
\lefteqn{
{\cal L}_{k}^{(n)} (\beta, s_n, \omega_n) \; = \; \frac
{ \Delta^k \, \{ (n + \beta)^{k-1} \> s_n / \omega_n\} }
{ \Delta^k \, \{ (n + \beta)^{k-1}  / \omega_n \} }
} \nonumber \\
& \; = \; \frac
{\displaystyle
\sum_{j=0}^{k} \; ( - 1)^{j} \; {{k} \choose {j}} \;
\frac
{(\beta + n +j )^{k-1}} {(\beta + n + k )^{k-1}} \;
\frac {s_{n+j}} {\omega_{n+j}} }
{\displaystyle
\sum_{j=0}^{k} \; ( - 1)^{j} \; {{k} \choose {j}} \;
\frac
{(\beta + n +j )^{k-1}} {(\beta + n + k )^{k-1}} \;
\frac {1} {\omega_{n+j}} } \, ,
\label{LevTr}
\\
\lefteqn{
{\cal S}_{k}^{(n)} (\beta, s_n, \omega_n) \; = \; \frac
{ \Delta^k \, \{ (n + \beta)_{k-1} \> s_n / \omega_n\} }
{ \Delta^k \, \{ (n + \beta)_{k-1}  / \omega_n \} }
} \nonumber \\
& \; = \; \frac
{\displaystyle
\sum_{j=0}^{k} \; ( - 1)^{j} \; {{k} \choose {j}} \;
\frac {(\beta + n +j )_{k-1}} {(\beta + n + k )_{k-1}} \;
\frac {s_{n+j}} {\omega_{n+j}} }
{\displaystyle
\sum_{j=0}^{k} \; ( - 1)^{j} \; {{k} \choose {j}} \;
\frac {(\beta + n +j )_{k-1}} {(\beta + n + k )_{k-1}} \;
\frac {1} {\omega_{n+j}} } \, .
\label{SidTr}
\end{eqnarray}
Here, $\{ s_n \}_{n=0}^{\infty}$ is the sequence to be transformed, and
$\{ \omega_n \}_{n=0}^{\infty}$ is a sequence of truncation error
estimates. The shift parameter $\beta$ has to be positive in order to
permit $n = 0$ in Eqs.\ (\ref{LevTr}) and (\ref{SidTr}). The most
obvious choice is $\beta = 1$. This has been the choice in virtually all
previous applications of these sequence transformations, and we will use
$\beta = 1$ unless explicitly stated. However, we show in Appendix
\ref{App_B} that in some cases it may be very advantageous to choose
other values for $\beta$.

The numerator and denominator sums in Eqs.\ (\ref{LevTr}) and
(\ref{SidTr}) can also be computed with the help of the three-term
recursions (Eq.\ (7.2-8) of Ref.\ \cite{We89})
\begin{eqnarray}
\lefteqn{L_{k+1}^{(n)} (\beta) \; = \; L_k^{(n+1)} (\beta)} \nonumber \\
& - & {\displaystyle \frac
{ (\beta + n + k ) (\beta + n + k)^{k-1} }
{ (\beta + n + k + 1)^k }} \, L_k^{(n)} (\beta)
\end{eqnarray}
and (Eq.\ (8.3-7) of Ref.\ \cite{We89})
\begin{eqnarray}
\lefteqn{S_{k+1}^{(n)} (\beta) \; = \; S_k^{(n+1)} (\beta)} \nonumber \\
& - & {\displaystyle \frac
{ (\beta + n + k ) (\beta + n + k - 1) }
{ (\beta + n + 2 k ) (\beta + n + 2 k - 1) }} \, S_k^{(n)} (\beta) \, .
\end{eqnarray}
The initial values $L_0^{(n)} (\beta) = S_0^{(n)} (\beta) = s_n /
\omega_n$ and $L_0^{(n)} (\beta) = S_0^{(n)} (\beta) = 1 / \omega_n$
produce the numerator and denominator sums, respectively, of ${\cal
L}_{k}^{(n)} (\beta , s_n, \omega_n)$ and ${\cal S}_{k}^{(n)} (\beta ,
s_n, \omega_n)$.

The performance of the sequence transformations ${\cal L}_{k}^{(n)}
(\beta, s_n, \omega_n)$ and ${\cal S}_{k}^{(n)} (\beta, s_n, \omega_n)$
depends crucially on the remainder estimates $\{ \omega_n
\}_{n=0}^{\infty}$. Under exceptionally favourable circumstances it may
be possible to construct explicit approximations to the truncation
errors of the input sequence $\{ s_n \}_{n=0}^{\infty}$. In most cases
of practical interest this is not possible, since only the numerical
values of a finite number of sequence elements is available. Hence, in
practice one has to determine the remainder estimates from these
numerical values.

On the basis of purely heuristic and asymptotic arguments, Levin
\cite{Le73} proposed some simple remainder estimates which -- when used
in Eq.\ (\ref{LevTr}) -- lead to Levin's $u$, $t$, and $v$
transformations. These remainder estimates can also be used in Eq.\
(\ref{SidTr}) (Section 8.4 of Ref.\ \cite{We89}).

However, in the case of a strictly alternating series the best {\em
simple\/} truncation error estimate is the first term neglected in the
partial sum (p.\ 259 of Ref.\ \cite{Kn64}). This is also the best simple
estimate for the truncation error of a strictly alternating
nonterminating hypergeometric series ${}_2 F_0 (\alpha, \beta; - z)$
with $\alpha, \beta, z > 0$, which diverges strongly for every $z
\ne 0$ (Theorem 5.12-5 of Ref.\ \cite{Ca77}). Consequently, in the case
of convergent or divergent alternating series it is a natural idea to
use the first term neglected in the partial sum as the remainder
estimate, as proposed by Smith and Ford \cite{SmFo79}.

In this article, we always accelerate the convergence of the partial
sums ${\bf S}_n$ of the Van Wijngaarden transformed series (\ref{A2B})
which is strictly alternating if the original series is monotone. Thus,
we use exclusively the sequence transformations (\ref{LevTr}) and
(\ref{SidTr}) in combination with the remainder estimate proposed by
Smith and Ford \cite{SmFo79}:
\begin{equation}
\omega_n \; = \; \Delta {\bf S}_n \; = \;
(-1)^{n+1} {\bf A}_{n+1} \, .
\label{d_Est}
\end{equation}
This yields the following sequence transformations (Eqs.\ (7.3-9) and
(8.4-4) of Ref.\ \cite{We89}):
\begin{eqnarray}
\lefteqn{d_{k}^{(n)} \bigl( \beta, {\bf S}_n \bigr) \; = \;
{\cal L}_{k}^{(n)} (\beta, {\bf S}_n, \Delta {\bf S}_n)} \nonumber \\
& \qquad = &
{\cal L}_{k}^{(n)} (\beta, {\bf S}_n, (-1)^{n+1} {\bf A}_{n+1}) \, ,
\label{dLevTr}
\\
\lefteqn{{\delta}_{k}^{(n)} \bigl( \beta, {\bf S}_n \bigr)
\; = \;
{\cal S}_{k}^{(n)} (\beta, {\bf S}_n, \Delta {\bf S}_n)} \nonumber \\
& \qquad = &
{\cal S}_{k}^{(n)} (\beta, {\bf S}_n, (-1)^{n+1} {\bf A}_{n+1}) \, .
\label{dSidTr}
\end{eqnarray}
Unless explicitly stated, we use $\beta = 1$. In the applications
described in the following, it will be clear from the context which
monotone series was transformed according to Van Wijngaarden such as to
produce the input data for the nonlinear sequence transformations which
are the partial sums ${\bf S}_n$ of the alternating series (\ref{A2B}).

Alternative remainder estimates for the sequence transformations
(\ref{LevTr}) and (\ref{SidTr}) were discussed in Sections 7 and 8 of
Ref.\ \cite{We89} or in Refs. \cite{We94a,HoWe95}.

>From a purely theoretical point of view, the sequence transformations
$d_{k}^{(n)}$ and ${\delta}_{k}^{(n)}$ as well as their parent
transformations ${\cal L}_{k}^{(n)}$ and ${\cal S}_{k}^{(n)}$ have the
disadvantage that no general convergence proof is known. Only for some
special model problems could rigorous convergence proofs be obtained
(Refs.\ \cite{Sid79,Sid80,Sid86} or Sections 13 and 14 of Ref.\
\cite{We89}). However, there is overwhelming empirical evidence that
$d_{k}^{(n)}$ and ${\delta}_{k}^{(n)}$ work very well in the case of
convergent or divergent alternating series for instance as they occur in
special function theory\cite{BreRZa91,We89,Le73,SmFo79,SmFo82,We90,
WeCi90,We96c,RoyBhaBho96,BhaBhoRoy97} or in quantum mechanical
perturbation theory
\cite{We90,WeCi90,WeCiVi91,We92,WeCiVi93,We96a,We96b,We96d,We70}.

Pelzl and King \cite{PeKi98} only used the Levin transformation
$d_{k}^{(n)} \bigl( \beta, {\bf S}_n \bigr)$ for the acceleration of the
convergence of Van Wijngaarden transformed alternating series. However,
we shall show that the closely related transformation
${\delta}_{k}^{(n)} \bigl( \beta, {\bf S}_n \bigr)$ frequently gives
better results.

Application of the sequence transformations (\ref{dLevTr}) and
(\ref{dSidTr}) to the partial sums ${\bf S}_n$ of the infinite series
(\ref{A2B}) produces doubly indexed sets of transforms, which depend on
the starting index $n$ of the input data and the transformation order
$k$:
\begin{equation}
\left\{ {\bf S}_n, {\bf S}_{n+1}, \ldots ,
{\bf S}_{n+k}, {\bf S}_{n+k+1} \right\} \to T_{k}^{(n)} \, .
\end{equation}
Here, $T_{k}^{(n)}$ stands for either $d_{k}^{(n)} \bigl( \beta, {\bf
S}_n \bigr)$ or ${\delta}_{k}^{(n)} \bigl( \beta, {\bf S}_n \bigr)$.

The transforms can be displayed in a rectangular scheme called the {\em
table\/} of the transformation. In principle, there is an unlimited
variety of different paths, on which we can move through the table in a
convergence acceleration or summation process. In this article, we
always proceed on such a path that for a given set of input data, we use
the transforms with the highest possible transformation order as
approximations to the limit of the series to be transformed. Thus, given
the input data $\left\{ {\bf S}_0, {\bf S}_{1}, \ldots , {\bf S}_{m},
{\bf S}_{m+1} \right\}$ we always use the sequence $\left\{ T_{0}^{(0)},
T_{1}^{(0)}, \ldots , T_{m}^{(0)} \right\}$ of transforms as
approximations to the limit (compare for instance Eqs.\ (7-5.4) and
(7-5.5) of Ref.\ \cite{We89}), where $T_{m}^{(0)}$ should provide the
best approximation since it has the highest transformation order.

In Appendix \ref{App_A}, we explain why the transformations ${\cal
L}_{k}^{(n)} (\beta, s_n, \omega_n)$ and ${\cal S}_{k}^{(n)} (\beta,
s_n, \omega_n)$ are much more effective in the case of an alternating
series than in the case of a monotone series.

\section{The Riemann zeta function}
\label{sec:ZetaFun}

The Riemann zeta function is discussed in most books on special
functions. Particularly detailed treatments can be found in monographs
by Edwards \cite{Edw74}, Titchmarsh \cite{Tit86}, Ivi\'{c} \cite{Ivi85},
and Patterson \cite{Pat88}. Many applications of the zeta function in
theoretical physics are described in books by Elizalde, Odintsov, Romeo,
Bytsenko, and Zerbini \cite{EliOdiRomBytZer94} and by Elizalde
\cite{Eli95}.

In this section, we discuss how the CNCT can be used for the evaluation
of the Riemann zeta function, starting from the logarithmically
convergent Dirichlet series (\ref{ZetaFun}). We do not claim that our
approach is necessarily more powerful than other, more specialized
techniques for the evaluation of the Riemann zeta function. However,
because of its simplicity, the Dirichlet series (\ref{ZetaFun}) is
particularly well suited for an illustration of our approach.

In the case of the Riemann zeta function, the terms ${\bf b}_{k}^{(j)}$
of the series (\ref{A2B}) for ${\bf A}_j$ are given by
\begin{equation}
{\bf b}_{k}^{(j)} \; = \; \frac {\left( 2^{1-z} \right)^k} {(j+1)^z} \, .
\label{zeta_b}
\end{equation}
Thus, the series (\ref{A2B}) is a power series in $2^{1-z}$ which
converges linearly. Moreover, it is essentially the geometric series in
the variable $2^{1-z}$ that can be expressed in closed form according to
\begin{equation}
{\bf A}_j \; = \;
\frac {1} {1-2^{1-z}} \, \frac {1} {(j+1)^z} \, .
\label{CondSumZeta}
\end{equation}
Inserting this into the infinite series on the right-hand side of Eq.\
(\ref{VWSerTran}) yields the following transformed alternating series
for the Riemann zeta function:
\begin{equation}
\zeta (z) \; = \; \frac {1} {1-2^{1-z}} \,
\sum_{j=0}^{\infty} \, \frac {(-1)^j} {(j+1)^z} \, .
\label{ZetaAlt}
\end{equation}
An alternative proof of the validity of this series expansion can be
found in Section 2.2 of Titchmarsh's book \cite{Tit86}.

The terms of the alternating series (\ref{ZetaAlt}) decay as slowly in
magnitude as the terms of the logarithmically convergent Dirichlet
series (\ref{ZetaFun}). Nevertheless, the series (\ref{ZetaAlt}) offers
some substantial computational advantages because it is alternating. For
example, it converges for ${\rm Re} (z) > 0$, whereas the Dirichlet
series converges only for ${\rm Re} (z) > 1$. For ${\rm Re} (z) \le 0$
both the Dirichlet series (\ref{ZetaFun}) and the alternating series
(\ref{ZetaAlt}) diverge. However, the alternating series can provide an
analytic continuation of the Riemann zeta function for ${\rm Re} (z) \le
0$ if it is combined with a suitable summation process. For this
purpose, we can use the same nonlinear transformations (\ref{dLevTr})
and (\ref{dSidTr}) as in the case of a {\em convergent\/} Van
Wijngaarden transformed alternating series.

This observation extends the applicability of the CNCT, which was
originally designed for convergent monotone series only. Of course, this
extension to divergent series is only possible if the monotone series
(\ref{A2B}) for ${\bf A}_j$ can be summed. Normally, the summation of a
monotone divergent series is very difficult. However, in the case of the
Riemann zeta function this is trivial since the series (\ref{A2B}) is --
as remarked above -- essentially the geometric series in the variable
$2^{1-z}$ which can be expressed in closed form according to Eq.\
(\ref{CondSumZeta}), no matter whether $\vert 2^{1-z} \vert < 1$ or
$\vert 2^{1-z} \vert > 1$ holds.

Our first numerical example is the zeta function with argument $z =
1.01$. The Dirichlet series (\ref{ZetaFun}) converges for this argument,
but its convergence is so slow that it is computationally useless. It
follows from the truncation error estimate (\ref{ZetaErrEst}) that on
the order of $10^{600}$ terms would be needed to achieve the modest
relative accuracy of $10^{-6}$ if the Dirichlet series were summed by
adding up the terms. A much more efficient evaluation of $\zeta(z)$ near
${\rm Re} (z) = 1$ is based on the Euler-Maclaurin formula (Eqs. (2.01)
and (2.02) on p.\ 285 of Ref.\ \cite{Olv74})
\begin{mathletters}
\label{EulerMaclaurin}
\begin{eqnarray}
& & \sum_{n=N}^{M} \, f (n) \; = \; \int_{N}^{M} \, f (x) \, {\rm d} x
\, + \, \frac{1}{2} \left[ f (N) + f (M) \right] \nonumber \\
& & \quad + \, \sum_{j=1}^{q} \, \frac {B_{2j}} {(2j)!}
\left[ f^{(2j-1)} (M) - f^{(2j-1)} (N) \right] \, + \, R_q \, ,
\\
& & R_q \; = \; - \frac{1}{(2q)!}
\int_{N}^{M} \, B_{2q}
\bigl( x - [\mkern - 2.5 mu  [x] \mkern - 2.5 mu ] \bigr) \,
f^{(2q)} (x) \, {\rm d} x \, .
\end{eqnarray}
\end{mathletters}
Here, $[\mkern - 2.5 mu  [x] \mkern - 2.5 mu ]$ is the integral part of
$x$, i.e., the largest integer $m$ satisfying $m \le x$, $B_k (x)$
is a Bernoulli polynomial defined by the generating function (Eq. (1.06)
on p.\ 281 of Ref.\ \cite{Olv74})
\begin{equation}
\frac {t \exp (x t)} {\exp (t) - 1} \; = \;
\sum_{m=0}^{\infty} \, B_m (x) \, \frac {t^m} {m!} \, ,
\qquad \vert t \vert < 2 \pi \, ,
\end{equation}
and
\begin{equation}
B_m \; = \; B_m (0)
\label{BernNum}
\end{equation}
is a Bernoulli number (p.\ 281 of Ref.\ \cite{Olv74}).

In the next step, we rewrite the Dirichlet series (\ref{ZetaFun}) as
follows,
\begin{eqnarray}
\zeta(z) \; = \; \sum_{m=0}^{N-1} \frac{1}{(m+1)^z} \, + \,
\sum_{m=N}^{\infty} \frac{1}{(m+1)^z} \, ,
\end{eqnarray}
and apply the Euler-Maclaurin formula (\ref{EulerMaclaurin}) to the
infinite series on the right-hand side. With $N = 100$, we easily obtain
the value of the zeta function accurate to more than 15 significant
decimal digits. The value of $\zeta(1.01)$ to 15 decimal places is
\begin{equation}
10^{-3}\,\zeta(1.01) \; = \; 0.100~577~943~338~497 \, .
\end{equation}
We present this result (as well as all other subsequent numerical
results) normalized to a number in the interval $(0,1)$ by a
multiplicative power of $10$.

In Table \ref{Tab_4_1}, we apply the Euler transformation
(\ref{EuTrans}) and the nonlinear sequence transformations
(\ref{dLevTr}) and (\ref{dSidTr}) to the partial sums of the alternating
series (\ref{ZetaAlt}) with $z = 1.01$. We list as a function of $n$ the
partial sums
\begin{equation}
{\bf S}_n \; = \; \frac {1}{1-2^{1-z}} \,
\sum_{k=0}^{n} \, \frac {(-1)^k}{(k+1)^z}
\end{equation}
of the alternating series (\ref{ZetaAlt}), the partial sums
\begin{equation}
{\bf E}_n \; = \; \frac {1}{1-2^{1-z}} \,
\sum_{k=0}^{n} \, \frac {1} {2^{k+1}} \,
\sum_{m=0}^{k} \, \frac {\displaystyle (-1)^m \, {{k} \choose {m}}}
{\displaystyle (m+1)^{z}}
\end{equation}
of the Euler transformed series (\ref{EuTrans}), and the nonlinear
sequence transformations $d_{n}^{(0)} \bigl( 1, {\bf S}_0 \bigr)$ and
${\delta}_{n}^{(0)} \bigl( 1, {\bf S}_0 \bigr)$.

\medskip
\begin{center}
{\bf \Large Insert Table \ref{Tab_4_1} here}
\end{center}
\medskip

The results in Table \ref{Tab_4_1} show that the Euler transformation is
in the case of the alternating series (\ref{ZetaAlt}) much less
effective than the nonlinear sequence transformations $d_{n}^{(0)}
\bigl( 1, {\bf S}_0 \bigr)$ and ${\delta}_{n}^{(0)} \bigl( 1, {\bf S}_0
\bigr)$. When we applied Wynn's epsilon algorithm \cite{Wy56a} to the
alternating series (\ref{ZetaAlt}), we also obtained clearly inferior
results. This is in agreement with the observations of Pelzl and King
\cite{PeKi98}.

Thus, in the following we will only consider the nonlinear transformations
(\ref{dLevTr}) and (\ref{dSidTr}) which use explicit remainder
estimates.

In Table \ref{Tab_4_1}, we present apparently redundant data since our
transformation results in the last two columns do not change for $n \ge
12$. However, we include these results here as well as below in order to
demonstrate that our transformation remains stable even if we increase
the transformation order.

If the argument $z$ is zero or a negative integer, the Riemann zeta
function satisfies (p.\ 807 of Ref.\ \cite{AbrSte72}):
\begin{mathletters}
\begin{eqnarray}
\lefteqn{} \nonumber \\
\zeta (0) & = & - \frac{1}{2} \, , \qquad \zeta (-2m) \; = \; 0 \, ,
\\
\zeta (1-2m) & = & - \, \frac {B_{2m}} {2m} \, ,
\qquad m \; = \; 1, 2, \ldots \, .
\label{BernoulZetaNegArg}
\end{eqnarray}
\end{mathletters}
Here, $B_{2m}$ is a Bernoulli number defined in Eq.\ (\ref{BernNum}).

Our second numerical example is the zeta function with argument $z = -
1$ which, because $B_2 = 1/6$, satisfies
\begin{equation}
\zeta (-1) \; = \; - 1/12 \; = \; -0.0833~333\ldots \, .
\end{equation}
As remarked above, the alternating series (\ref{ZetaAlt}) diverges for
$z = - 1$. However, our results in Table \ref{Tab_4_2} show that this
series can be summed effectively by the nonlinear sequence
transformations (\ref{dLevTr}) and (\ref{dSidTr}).

\medskip
\begin{center}
{\bf \Large Insert Table \ref{Tab_4_2} here}
\end{center}
\medskip

The results in Table \ref{Tab_4_2} indicate that the sequence
transformation ${\delta}_{n}^{(0)} \bigl(1, {\bf S}_0 \bigr)$ is exact
for the alternating series (\ref{ZetaAlt}) with $z = -1$. The exactness
of the sequence transformations ${\cal L}_{k}^{(n)} (\beta, s_n,
\omega_n)$, Eq.\ (\ref{LevTr}), and ${\cal S}_{k}^{(n)} (\beta, s_n,
\omega_n)$, Eq.\ (\ref{SidTr}), for the infinite series (\ref{ZetaAlt})
and for related series is discussed in Appendix \ref{App_B}.

There is considerable research going on in connection with the
Riemann-Siegel conjecture \cite{BerKea92,Par94,Ber95,BhaKhaReiTom97}. In
this context, it is necessary to evaluate the Riemann zeta function on
the so-called critical line $z = 1/2 + t~i$ ($- \infty < t < \infty$).
This can also be accomplished by applying the sequence transformations
(\ref{dLevTr}) and (\ref{dSidTr}) to the partial sums of the alternating
series (\ref{ZetaAlt}). In Table \ref{Tab_4_3}, where we treat the case
$z = 1/2 + 13.7~i$, the Levin transformation outperforms the Weniger
transformation. For limitations of space, we do not present the partial
sums of the transformed alternating series in Table
\ref{Tab_4_3}.

\medskip
\begin{center}
{\bf \Large Insert Table \ref{Tab_4_3} here}
\end{center}
\medskip

\section{The Lerch Transcendent and Related Functions}
\label{sec:LerTran}

The Lerch transcendent (p.\ 32 of Ref.\ \cite{MaObSo66})
\begin{equation}
\Phi (z, s, \alpha) \; = \; \sum_{n=0}^{\infty} \,
\frac {z^n} {(\alpha + n)^s} \, , \qquad \vert z \vert < 1 \, ,
\label{Lerch}
\end{equation}
contains many special functions as special cases, for example the
Riemann zeta function
\begin{equation}
\zeta (s) \; = \; \Phi (1, s, 1)
\end{equation}
or the Jonqui\`{e}re function (p.\ 33 of Ref.\ \cite{MaObSo66})
\begin{equation}
F (z, s) \; = \; \sum_{n=0}^{\infty} \, \frac {z^{n+1}} {(n+1)^s}
\; = \; z \, \Phi (z, s, 1) \, .
\end{equation}
In the physics literature, the Jonqui\`{e}re function is usually called
generalized logarithm or polylogarithm, and the following notation is
used:
\begin{equation}
{\rm Li}_s (z) \; = \; \sum_{k=0}^{\infty} \, \frac{z^{k+1}}{(k+1)^s}
\; = \; F(z, s) \, .
\label{PolyLog}
\end{equation}
This is also the notation and terminology used in the book by Lewin
\cite{Lew81}.

Lerch transcendents and their special cases, the polylogarithms, are
ubiquitous in theoretical physics. They play a role in Bose-Einstein
condensation \cite{GriWuStr97,MacL97}, and they are particularly
important in quantum field theory. For example, they occur in integrands of
one-dimensional numerical integrations in quantum electrodynamic bound
state calculations \cite{Pa92,Pa93,JePa96,JeSoMo97,IndMo98}.

The series expansion (\ref{Lerch}) for $\Phi (z, s, \alpha)$ converges
very slowly if the argument $z$ is only slightly smaller than one.
However, the CNCT can be used for an efficient and reliable evaluation
of the Lerch transcendent and its special cases.

In our first example (Table \ref{Tab_5_1}), we evaluate
\begin{eqnarray}
\label{example:li1}
& & 10^{-2}\,{\rm Li}_1 (0.99999) =
- 10^{-2}\,\ln(0.00001) \nonumber\\
& & \;\;\;\;\;\; = 0.115~129~254~649~702 \, .
\end{eqnarray}
The series (\ref{PolyLog}) for ${\rm Li}_s (z)$ converges linearly for
$\vert z \vert < 1$, but
\begin{equation}
{\rm Li}_1 (z) \; = \; \sum_{k=0}^{\infty} \frac{z^{k+1}}{k+1}
\; = \; - \ln(1 - z)
\label{SerLi_1}
\end{equation}
has a singularity at $z = 1$. Thus, in Table \ref{Tab_5_1} we evaluate
${\rm Li}_1 (z)$ in the immediate vicinity of a singularity.

\medskip
\begin{center}
{\bf \Large Insert Table \ref{Tab_5_1} here}
\end{center}
\medskip

In contrast, the series expansions
\begin{equation}
{\rm Li}_2 (z) \; = \; \sum_{k=0}^{\infty} \frac{z^{k+1}}{(k+1)^2}
\label{SerLi_2}
\end{equation}
and
\begin{equation}
{\rm Li}_3 (z) \; = \; \sum_{k=0}^{\infty} \frac{z^{k+1}}{(k+1)^3}
\label{SerLi_3}
\end{equation}
converge logarithmically for $z = 1$. In Tables \ref{Tab_5_2} and
\ref{Tab_5_3}, we consider
\begin{equation}
10^{-1}\,{\rm Li}_2(0.99999) \; = \; 0.164~480~893~699~293\dots
\end{equation}
and
\begin{equation}
10^{-1}\,{\rm Li}_3(0.99999) = 0.120~204~045~438~733 \, ,
\end{equation}
respectively. As in most previous examples, the Weniger transform
${\delta}_{n}^{(0)} \bigl(1, {\bf S}_0 \bigr)$ produces in Tables
\ref{Tab_5_2} and \ref{Tab_5_3} somewhat better results than the Levin
transform $d_{n}^{(0)} \bigl(1, {\bf S}_0 \bigr)$.

\medskip
\begin{center}
{\bf \Large Insert Tables \ref{Tab_5_2} and \ref{Tab_5_3} here}
\end{center}
\medskip

In our last example of this section, we evaluate
\begin{equation}
10^4\,\Phi(0.99999,2,10000) \; = \; 0.798~585~139~222~548
\end{equation}
with the help of the CNCT. The prefactor of $10000$ is introduced so
that the final result is of order one.

\medskip
\begin{center}
{\bf \Large Insert Table \ref{Tab_5_4} here}
\end{center}
\medskip

\section{The Generalized Hypergeometric Series}
\label{sec:GenHygFun}

The Gaussian hypergeometric function ${}_2 F_1 (a, b; c; z)$, defined by
the series expansion (\ref{2F1_Ser}) for $\vert z \vert < 1$, is one of
the most used and best understood special functions. It is discussed in
virtually all books on special functions. A particularly detailed
treatment can be found in the book by Slater \cite{Sla66}.

The fact that the mathematical properties of the hypergeometric function
${}_2 F_1$ are so well understood, greatly facilitates its computation.
The series (\ref{2F1_Ser}) does not suffice for the computation of a
nonterminating ${}_2 F_1$ since it converges only if $\vert z \vert <
1$. By contrast, the hypergeometric function ${}_2 F_1$ is a multivalued
function defined in the whole complex plane with branch points at $z =
1$ and $\infty$.

However, a hypergeometric function ${}_2 F_1 (a, b; c; z)$ can be
transformed into the sum of two other ${}_2 F_1$'s with arguments $w = 1
- z$, $w = 1/z$, $w = 1/(1-z)$, or $w = 1-1/z$, respectively. Thus, the
argument $w$ of the two resulting hypergeometric series can normally be
chosen in such a way that the two new series in $w$ either converge, if
the original series in $z$ diverges, or they converge more rapidly if
the original series converges too slowly to be numerically useful.

For example, if $\vert 1-z \vert < 1$ and if $c - a - b$ is not a
positive or negative integer, then we can use the analytic continuation
formula (Eq.\ (15.3.6) of Ref.\ \cite{AbrSte72})
\begin{eqnarray}
\lefteqn{{}_2 F_1 (a, b; c; z)} \nonumber \\
& = &
\frac {\Gamma (c) \Gamma (c-a-b)} {\Gamma (c-a) \Gamma (c-b)} \,
{}_2 F_1 (a, b; a+b-c+1; 1-z) \nonumber \\
& & + \,
\frac {\Gamma (c) \Gamma (a+b-c)} {\Gamma (a) \Gamma (b)} \,
(1-z)^{c-a-b} \nonumber \\
& & \times \, {}_2 F_1 (c-a, c-b; c-a-b+1; 1-z) \, .
\label{AC_1}
\end{eqnarray}
Obviously, the two hypergeometric series in $1-z$ will converge rapidly
in the vicinity of $z = 1$. With the help of this or other analytic
continuation formulas it is possible to compute a hypergeometric
function ${}_2 F_1 (a, b; c; z)$ with arbitrary real argument $z \in (-
\infty, + \infty)$ effectively via the resulting hypergeometric series
(see p.\ 127 of Ref.\ \cite{Tem96} or Table I of Ref.\ \cite{For97}).

The situation is much less favourable in the case of the generalized
hypergeometric series
\begin{equation}
{}_3 F_2 \bigl( a, b, c; d, e; z \bigr) \; = \;
\sum_{m=0}^{\infty} \, \frac{(a)_m (b)_m (c)_m}{(d)_m (e)_m} \,
\frac{z^m}{m!} \, ,
\label{Def_3F2}
\end{equation}
and it gets progressively worse in the case of higher generalized
hypergeometric series ${}_{p+1} F_p$ with $p \ge 3$. The analytic
continuation formulas of the type of Eq.\ (\ref{AC_1}) are not always
known, and those that are known become more and more complicated with
increasing $p$ \cite{Bue87a,Bue87b,Bue88,Bue92}. Thus, the efficient and
reliable computation of a generalized hypergeometric function ${}_{p+1}
F_p$ with $p \ge 2$ via its series expansion (\ref{GenHygSer}) is still
a more or less open problem. For example, in Theorem 2 of Ref.\
\cite{Bue92} explicit expressions for the analytic continuation formulas
of a ${}_{p+1} F_p$ with $p = 2, 3, \ldots$ around $z = 1$ were
constructed. However, these expansions are by some orders of magnitude
more complicated than the analogous formula (\ref{AC_1}) for a 
${}_2 F_1$.

We show here that the CNCT can be very useful if the argument $z$ of a
generalized hypergeometric series ${}_{p+1} F_p$ is only slightly
smaller than one, even if $z = 1$ is a singularity.

The function ${}_3 F_2(1,3/2,5;9/8,47/8;z)$ has a singularity at $z=1$
because the sum of the real parts of the numerator parameters is greater
than the sum of the real parts of the denominator parameters (see p.\ 45
of Ref.\ \cite{Sla66}). Thus, in Table \ref{Tab_6_1} we evaluate the
function
\begin{eqnarray}
& & 10^{-4} \, {}_3 F_2(1,3/2,5;9/8,47/8;0.99999)
\nonumber\\
& & \;\;\;\;\;\; = 0.238~434~298~763~330
\end{eqnarray}
in the immediate vicinity of the singularity. With the help of the CNCT,
it is nevertheless possible to evaluate this function without any
noticeable loss of significant digits.

\medskip
\begin{center}
{\bf \Large Insert Table \ref{Tab_6_1} here}
\end{center}
\medskip

The function ${}_3 F_2(1,3,7;5/2,14;z)$ does not have a singularity
at $z=1$, i.e., the hypergeometric series (\ref{Def_3F2}) converges for
$z = 1$, albeit very slowly. In Table \ref{Tab_6_2} we compute the
function
\begin{eqnarray}
& & 10^{-1}\,{}_3 F_2(1,3,7;5/2,14;0.99999)
\nonumber\\
& & \;\;\;\;\;\; = 0.267~102~823~984~762
\end{eqnarray}
with the help of the CNCT.

\medskip
\begin{center}
{\bf \Large Insert Table \ref{Tab_6_2} here}
\end{center}
\medskip

There is another major difference between a Gaussian hypergeometric
series ${}_2 F_1$ and a generalized hypergeometric series ${}_{p+1} F_p$
with $p \ge 2$. The series (\ref{2F1_Ser}) for a hypergeometric function
${}_2 F_1 (a, b; c; z)$ with unit argument $z = 1$ converges if ${\rm
Re} (c) > {\rm Re}(a + b)$, and its value is given by the Gauss
summation theorem (compare for instance Section 1.7 of Ref.\
\cite{Sla66}):
\begin{equation}
{}_2 F_1 (a, b; c; 1) \; = \;
\frac {\Gamma (c) \Gamma (c-a-b)} {\Gamma (c-a) \Gamma (c-b)} \, .
\label{GaussSumTheo}
\end{equation}
The series (\ref{GenHygSer}) for ${}_{p+1} F_p \bigl(\alpha_1, \ldots,
\alpha_{p+1}; \beta_1, \ldots, \beta_p; z \bigr)$ with $p \ge 2$
converges at unit argument $z = 1$ if ${\rm Re} (\beta_1 + \ldots +
\beta_p) > {\rm Re} (\alpha_1 + \ldots + \alpha_{p+1})$ (compare
p.\ 45 of Ref.\ \cite{Sla66}). In Theorem 1 of Ref.\ \cite{Bue92}, an
analogue of the Gauss summation theorem for a generalized hypergeometric
series ${}_{p+1} F_p$ with unit argument was formulated. However, this
expression is not a simple ratio of gamma functions as in Eq.\
(\ref{GaussSumTheo}), but an infinite series with complicated terms.

Otherwise, there is a variety of different {\em simple\/} summation
theorems for generalized hypergeometric series ${}_{p+1} F_p$ with unit
argument which are ratios of gamma function. However, these summation
theorems depend upon the value of $p$, and are usually only valid for
certain values and/or combinations of the parameters $\alpha_1$,
$\alpha_2$, \ldots, $\alpha_{p+1}$ and $\beta_1$, $\beta_2$, \ldots,
$\beta_p$. A list of these simple summation theorems can be found in
Appendix III of Slater's book \cite{Sla66}. The reconstruction of known
simple summation theorems for nonterminating generalized hypergeometric
series and the construction of new ones with the help of computer
algebra systems is discussed in books by Petkov\v{s}ek, Wilf, and
Zeilberger \cite{PetWilZei96} and Koepf \cite{Koe98}.

With the help of Watson's summation theorem (p.~245 Ref. \cite{Sla66})
we obtain
\begin{eqnarray}
& & 10^{-1}\,{}_3 F_2(1,3,7;5/2,14;1) =
  \frac{567567 \, \pi^2}{20971520}
\nonumber\\
& & \;\;\;\;\;\; = 0.267~108~047~538~428 \, .
\end{eqnarray}
We show in Table \ref{Tab_6_3} that this result can also be obtained
with the help of the CNCT, starting from the hypergeometric series
(\ref{Def_3F2}).

\medskip
\begin{center}
{\bf \Large Insert Table \ref{Tab_6_3} here}
\end{center}
\medskip

\section{Products of Bessel and Hankel Functions}
\label{sec:ProdBessHank}

In this section, we want to apply the CNCT to series whose terms involve
more complicated entities and which are more typical of the problems
encountered in theoretical physics.

As an example we consider here a series whose terms are products of
spherical Bessel and Hankel functions (p.\ 435 of Ref.\ 
\cite{AbrSte72}):
\begin{eqnarray}
j_l (z) & = & [\pi/(2 z)]^{1/2} \, J_{l+1/2} (z) \, ,
\\
h^{(1)}_l (z) & = & [\pi/(2 z)]^{1/2} \, H^{(1)}_{l+1/2} (z) \, .
\end{eqnarray}
Here, the index $l$ is a nonnegative integer, and $J_{l+1/2}$ and
$H^{(1)}_{l+1/2}$ are Bessel and Hankel functions, respectively (pp.\
358 and 360 of Ref.\ \cite{AbrSte72}). In physical applications, the
index $l$ usually finds a natural interpretation as an angular momentum
quantum number.

Here, we investigate the following model problem:
\begin{equation}
\label{ht1}
\frac{\exp\bigl( -y[1-r] \bigr)}{y[1-r]} \; = \; - \,
\sum_{l=0}^{\infty} (2l+1) \, j_l (iry) \, h^{(1)}_l (iy) \, ,
\end{equation}
where $y$ is positive and $0 < r < 1$. The spherical Bessel and Hankel
functions $j_l (iry)$ and $h^{(1)}_l (iy)$ can also be expressed in
terms of {\em modified\/} spherical Bessel functions. However, we prefer
to retain the given notation because there are some inconsistencies in
the literature regarding the prefactors to be assigned to the {\em
modified} spherical Bessel functions, whereas the spherical Bessel and
Hankel functions are defined consistently in most textbooks.

With the help of known properties of Bessel functions and Legendre
polynomials it can be shown easily that the series expansion (\ref{ht1})
is a special case of the well known addition theorem of the Yukawa
potential (p.\ 107 of Ref.\ \cite{MaObSo66}),
\begin{eqnarray}
\lefteqn{\frac {\exp (-\gamma w)} {w} \; = \; (r \rho)^{-1/2}}
\nonumber \\
& & \quad \times \, \sum_{l=0}^{\infty} \, (2l+1) \, P_l (\cos \phi) \,
I_{l+1/2} (\gamma x) \, K_{l+1/2} (\gamma y) \, ,
\end{eqnarray}
where $w = (x^2 + y^2 - 2 x y \cos \phi)^{1/2}$, $0 < x < y$ and $\gamma
> 0$. This addition theorem is also the Green's function of the
three-dimensional modified Helmholtz equation.

First, we want to analyze the behaviour of the terms of the series on
the right-hand side of Eq.\ (\ref{ht1}) if the index $l$ becomes large.
The leading orders of the asymptotic expansions of the spherical Bessel
and Hankel functions $j_l (iry)$ and $h^{(1)}_l (iy)$ as $l \to \infty$
with $r$ and $y$ fixed can be obtained easily from the series expansion
for $j_l$ (Eq.\ (10.1.2) of Ref.\ \cite{AbrSte72} or Eq.\ (E.17) of
Ref.\ \cite{Mo74:2}) and from the explicit expression for $h^{(1)}_l$
(Eqs.\ (10.1.3) and (10.1.16) of Ref.\ \cite{AbrSte72} or Eq.\ (E.18) of
Ref.\ \cite{Mo74:2}), yielding
\begin{eqnarray}
j_l(iry) & \sim & \frac{(iry)^l}{(2 l + 1)!!} \,
\left[1 + {\rm O}\left( 1/l \right) \right] \, ,
\\
h^{(1)}_l(iy) & \sim & -i\,\frac{(2 l - 1)!! {\rm e}^{-y}}{(iy)^{l + 1}} 
\, \left[1 + {\rm O} \left( 1/l \right) \right] \, .
\end{eqnarray}
Thus, we obtain for the leading order of the asymptotic expansion ($l
\to \infty$) of the product of the spherical Bessel and the spherical
Hankel function:
\begin{equation}
\label{leading:order}
j_l (iry) \, h^{(1)}_l (iy) \; \sim \; - 
\frac{r^l {\rm e}^{-y}}{(2l + 1)y} \,
\left[1 + {\rm O}\left( 1/l \right) \right] \, .
\end{equation}
This asymptotic result shows that the infinite series on the right-hand
side of Eq.\ (\ref{ht1}) converges only for $\vert r \vert < 1$, and
that the function $\exp ( -y[1-r])/(y[1-r])$ has a singularity at $r =
1$.

Thus, it is interesting to use the CNCT for an evaluation of the series
on the right-hand side of Eq.\ (\ref{ht1}) in the immediate vicinity of
the singularity. This series is similar to the sum over angular momenta
encountered in QED bound state calculations (cf.\ Ref.\
\cite{JeSoMo98} and Eq.\ (4.2) in Ref.\ \cite{Mo74:2}).

In principle, the spherical Bessel and Hankel functions in the products
$j_l(iry) \, h^{(1)}_l(iy)$ can be evaluated recursively. If, however,
the CNCT is used, this is not useful. The products belonging to {\em
all} angular momenta $l$ are not needed. Instead, only the products for
some {\em specific} angular momenta have to be evaluated. Accordingly,
we prefer to evaluate the spherical Bessel function $j_l(iry)$ from its
series expansion (Eq.\ (10.1.2) of Ref.\ \cite{AbrSte72} or Eq.\ (E.17)
of Ref.\ \cite{Mo74:2}) and the spherical Hankel function
$h^{(1)}_l(iy)$ via its explicit expression (Eq.\ (10.1.16) of Ref.\
\cite{AbrSte72} or Eq.\ (E.18) of Ref.\ \cite{Mo74:2}).  

It should be taken into account that the recursive evaluation of the
terms of angular momentum decompositions is normally possible only in
the case of simple model problems. We shortly outline the applications
of the CNCT to a more involved computational problem occurring in QED
bound state calculations \cite{JeSoMo98}. The series to be evaluated in
the QED calculation is given by Eq.\ (4.2) in Ref.\ \cite{Mo74:2}. This
series can be rewritten as
\begin{equation}
\label{ang:sum}
S(r,y,t,\gamma) \; = \; J(r,y,t,\gamma) \, 
\sum_{\kappa=1}^{\infty} T_{\kappa}(r,y,t,\gamma) \, ,
\end{equation}
where $J(r,y,t,\gamma)$ has a simple mathematical structure. The term
$T_{\kappa}(r,y,t,\gamma)$ ($\kappa = 1,2,\dots$) can be rewritten as
(see Eq.\ (5.7) in Ref.\ \cite{Mo74:1})
\begin{eqnarray}
\label{Tk}
\lefteqn{T_{\kappa}(r,y,t,\gamma) \; = \;
\sum_{i,j=1}^{2} f_i\left(\frac{r\,y}{a}\right) \,
G_{\rm B,\kappa}^{ij}
\left[\frac{r\,y}{a},\frac{y}{a},\frac{i}{2}\,
\left(\frac{1}{t} - t\right)\right]} \nonumber\\
&& \times \,
f_j\left(\frac{y}{a}\right) \, 
A_{\kappa}\left(\frac{r\,y}{a},\frac{y}{a}\right) 
- f_{3-i}\left(\frac{r\,y}{a}\right) \nonumber\\ 
&& \times \,                   
G_{\rm B,\kappa}^{ij}\left[\frac{r\,y}{a},\frac{y}{a},\frac{i}{2}\,
\left(\frac{1}{t} - t\right)\right] \,
f_{3-j}\left(\frac{y}{a}\right) \, 
A^{ij}_{\kappa}\left(\frac{r\,y}{a},\frac{y}{a}\right) \, ,
\end{eqnarray}
where
\begin{itemize}
\item the functions $f_i$ are the radial Dirac-Coulomb wave functions
(see Eq.\ (A.8) in Ref.\ \cite{Mo74:1}),
\item the functions $G^{ij}_{\rm B,\kappa}$ are related to the radial
components of the relativistic Green function of the bound electron
(see Eq.\ (5.5) in Ref.\ \cite{Mo74:1}), and
\item the functions $A_{\kappa}$ and $A^{ij}_{\kappa}$ are related to
the angular momentum decomposition of the Green function of the
virtual photon and are defined in Eq.\ (5.8) of Ref.\ \cite{Mo74:1}.
\end{itemize} 
The quantity $a$ is a scaling variable for the subsequent radial
integration over $y$ (see Eq.\ (4.1) in Ref.\ \cite{Mo74:2}). The
variable $r$ denotes the ratio of the two radial coordinates ($0 < r <
1$), and $t$ is related to the (complex) energy of the virtual photon
($0 < t < 1$). The dependence on the coupling $\gamma$ of the electron
to the central field, which appears on the left-hand side of Eq.\
(\ref{Tk}), is implicitly contained in the functions $G^{ij}_{\rm
B,\kappa}$ on the right-hand side of Eq.\ (\ref{Tk}) (see Eq.\ (A.16) in
Ref.\ \cite{Mo74:1}). The propagator $G_{\rm B}$ in Eq.\ (\ref{Tk})
contains the relativistic Dirac-Coulomb Green function, whose radial
components are given in Eq.\ (A.16) in Ref.\ \cite{Mo74:1}. Recurrence
formulas relating the $G_{\rm B,\kappa}^{ij}$ for different values of
$\kappa$ are not known. Therefore the terms in the series
(\ref{ang:sum}) cannot be evaluated by recursion, but each one of these
terms can be computed using techniques described in Refs.\
\cite{JeSoMo98,Mo74:2}. Note that the angular momentum decomposition of
Eq.\ (\ref{ang:sum}) does not correspond to an expansion in terms of the
QED perturbation theory. The perturbative parameter in QED is the
elementary charge $e$ or the fine structure constant\footnote{We use
natural units as it is customary for QED bound state calculations
($\hbar = c = \epsilon_0 = 1$).}  $\alpha = e^2/(4\,\pi)$. The series
(\ref{ang:sum}) occurs in the evaluation of an energy shift (the Lamb
shift) of atomic levels due to the self energy of the electron. The self
energy is an effect described by second-order perturbation theory within
the framework of QED. The asymptotic behaviour of the terms $T_{\kappa}$
in Eq.\ (\ref{ang:sum}),
\begin{equation}
T_{\kappa}\; \sim \;  \frac{r^{2\,\kappa}}{\kappa} \,
\left[1 + {\rm O}\left(\frac{1}{\kappa}\right)\right] \, ,
\qquad \kappa \to \infty \, ,
\end{equation}
is similar to that of the model series (\ref{leading:order}) and to that
of the ${\rm Li}_1$-series discussed in Sec.\ \ref{sec:LerTran} (see
Eq.\ (\ref{example:li1}) and the results presented in Table 
\ref{Tab_5_1}). Slow convergence of the series in Eq.\ (\ref{ang:sum})
is observed in the limit $r \to 1$.

Returning to our model problem given in Eq.\ (\ref{ht1}), we consider
here the example of $r = 0.9999$ and $y = 0.7$. This yields
\begin{eqnarray}
& & 10^{-5}\,\sum_{l=0}^{\infty} (2\,l+1)\,j_l(i\,0.9999\times0.7)\,
h^{(1)}_l(i\,0.7) \nonumber\\
& & \;\;\;\;\; = \; - 0.142~847~143~207~135 \, .
\end{eqnarray}
Evaluation of the series over the spherical Bessel and Hankel functions
via the CNCT is presented in Table \ref{Tab_7_1}. If this series is
evaluated by adding up its terms, then about $450~000$ products of
Bessel and Hankel functions have to be evaluated for a relative accuracy
of $10^{-16}$ in the final result. This is to be compared with the
approximately $300$ evaluations of products of Bessel and Hankel
functions which are needed to compute the 25 Van Wijngaarden terms
presented in Table \ref{Tab_7_1} with the necessary accuracy. According
to our experience, the use of the CNCT in QED bound state calculations
reduces the amount of computer time for slowly convergent series of the 
type of Eq.\ (\ref{ang:sum}) by three orders of magnitude.

\medskip
\begin{center}
{\bf \Large Insert Table \ref{Tab_7_1} here}
\end{center}
\medskip

\section{Summary and Conclusions}
\label{sec:SumConclu}

Our approach for the acceleration of the convergence of monotone series
consists of two steps. In the first step, the Van Wijngaarden
condensation transformation, Eqs. (\ref{VWSerTran}) - (\ref{B2a}), is
applied to the original monotone series. In the second step, the
convergence of the resulting alternating series (\ref{A2B}) is
accelerated with the help of suitable nonlinear sequence
transformations.

Daniel \cite{Dan69} showed that the transformed alternating series
(\ref{A2B}) converges under relatively mild conditions to the same limit
as the original series. The condensation transformation cannot solve the
problems due to slow convergence on its own since the transformed
alternating series does not converge more rapidly than the original
series (see Eq.\ (\ref{aBoundA})). However, the convergence of
alternating series can be accelerated much more effectively than the
convergence of the original monotone series. Moreover, the
transformation of alternating series by sequence transformations is
numerically a relatively stable process, whereas the direct
transformation of monotone series inevitably leads to the loss of
accuracy due to round-off.

Many sequence transformations are known which are in principle capable
of accelerating the convergence of alternating series
\cite{Wi81,BreRZa91,We89}. However, their efficiency varies
considerably. In all numerical examples studied in the context of this
paper, we found that certain sequence transformations \cite{We89,Le73}
(which also use explicit estimates for the truncation errors as input
data) are significantly more powerful than other, better known sequence
transformations as for instance Wynn's epsilon algorithm \cite{Wy56a} or
the Euler transformation \cite{Euler1755} (which only use the partial
sums of the infinite series to be transformed as input data).
Consequently, we use in this article only the sequence transformations
${\cal L}_{k}^{(n)} (\beta, s_n, \omega_n)$ and ${\cal S}_{k}^{(n)}
(\beta, s_n, \omega_n)$ defined in Eqs.\ (\ref{LevTr}) and (\ref{SidTr})
in combination with the truncation error estimate (\ref{d_Est}). Those
properties of these sequence transformations, which are needed to apply
them successfully in a convergence acceleration process, are described
in Section \ref{sec:NonSeqTrans}. In Appendix \ref{App_A}, we explain
why these transformations are much more effective in the case of
alternating series than in the case of monotone series, and in Appendix
\ref{App_B} we discuss exactness properties of these transformations.

We consistently observed that the combined transformation is a
remarkably stable numerical process and entails virtually no loss of
numerical significance at intermediate stages. The evaluation of the
condensed series (\ref{B2a}) is a computationally simple task, provided
that the evaluation of terms of sufficiently high indices in the
original series is feasible. The nonlinear sequence transformations used
in the second step are also remarkably stable. The final results, which
we present in the tables, indicate that a relative accuracy of
$10^{-14}$ can be achieved with floating point arithmetic of 16 decimal
digits. Moreover, our method remains stable in higher transformation
orders.

In Sections \ref{sec:ZetaFun} - \ref{sec:GenHygFun} we demonstrate that
our approach is very useful for the computation of special functions
which are defined by logarithmically convergent series. Our approach
also works very well if the argument of a power series with monotone
coefficients is very close to the boundary of the circle of convergence
(even in the vicinity of singularities).

In Section \ref{sec:ZetaFun}, we treat the Dirichlet series
(\ref{ZetaFun}) for the Riemann zeta function which because of its
simplicity is excellently suited to illustrate the mechanism of our
combined transformation. The application of the Van Wijngaarden
transformation to the logarithmically convergent series (\ref{ZetaFun})
yields the known alternating series (\ref{ZetaAlt}) which does not
converge more rapidly than the original monotone series (\ref{ZetaFun}).
However, the alternating series can be transformed effectively by the
sequence transformations (\ref{dLevTr}) or (\ref{dSidTr}) in the case of
slow convergence or even summed in the case of divergence. Consequently,
the Riemann zeta function can be evaluated effectively and reliably even
if its argument $z$ is only slightly larger than one (Table
\ref{Tab_4_1}), if both the monotone series (\ref{ZetaFun}) and the
alternating series (\ref{ZetaAlt}) diverge (Table \ref{Tab_4_2}), or if 
the argument $z$ is complex (Table \ref{Tab_4_3}).

In Section \ref{sec:LerTran} we treat the Lerch transcendent $\Phi (z,
s, \alpha)$ which contains many other special functions as special
cases, for example the Riemann zeta function or the polylogarithms. In
Tables \ref{Tab_5_1} - \ref{Tab_5_4} we show that the combined
transformation makes it possible to evaluate $\Phi (z, s, \alpha)$
effectively and reliable via the power series (\ref{Lerch}) even if $z$
is very close to the boundary of the circle of convergence.

In Section \ref{sec:GenHygFun} we discuss the evaluation of a
nonterminating generalized hypergeometric series ${}_{p+1} F_p$ via its
defining power series (\ref{GenHygSer}) which converges for $\vert z
\vert < 1$. In the case of a Gaussian hypergeometric series ${}_2 F_1$,
explicit analytic continuation formulas are known which transform a
hypergeometric series with argument $z$ into the sum of two
hypergeometric series ${}_2 F_1$ with argument $1-z$. Thus, if the
convergence of a hypergeometric series ${}_2 F_1$ is slow because its
argument $z$ is only slightly smaller then one, then the two transformed
hypergeometric series ${}_2 F_1$ with argument $1-z$ will converge
rapidly. Moreover, if a Gaussian hypergeometric series ${}_2 F_1$ with
unit argument $z = 1$ converges, then its value is given by the Gauss
summation theorem (\ref{GaussSumTheo}). In the case of a generalized
hypergeometric series ${}_{p+1} F_p$ with $p \ge 2$, the situation is
much more difficult. Explicit analytic continuation formulas are either
unknown or they become increasingly complicated with increasing
$p$. Similarly, simple analogues for the Gauss summation theorem
(\ref{GaussSumTheo}) exist only for certain values of $p$ and for
certain combinations of the parameters of the ${}_{p+1} F_p$. In Tables
\ref{Tab_6_1} - \ref{Tab_6_3} we show that a generalized hypergeometric
series ${}_{p+1} F_p$ ($p \ge 2$) with an argument $z$ that is only
slightly smaller than one or with unit argument ($z = 1$) can be computed
effectively and reliably with the help of the combined transformation.

Partial wave decompositions of Green's functions, which occur for
example in quantum electrodynamic bound state calculations, entail more
complex mathematical entities than series expansions for special
functions. Nevertheless, we show in Section \ref{sec:ProdBessHank} that
the combined transformation can be applied successfully to these partial
wave decompositions. The series over products of spherical Bessel and
Hankel functions considered in this paper serves as a model problem for
the angular momentum decomposition of more complex Green's functions, as
for example the relativistic Green's function of the bound electron
\cite{Mo74:1}. In the model problem studied here as well as in a recent
evaluation of self energy corrections in bound systems \cite{JeSoMo98},
we observed a reduction in computer time by three orders of
magnitude. Thus, the combined transformations makes extensive and highly
accurate calculations feasible in situations that could otherwise be too
time-consuming.

\acknowledgments

We thank Prof.\ F.\ W.\ J.\ Olver for illuminating discussions. U.D.J.\
thanks J.~Baker, J.~Devaney and J.~Sims for stimulating discussions, and
acknowledges support from the Deutsche Forschungsgemeinschaft (contract
no.\ SO333-1/2) and the DAAD. P.J.M.\ acknowledges continued support by
the Alexander von Humboldt Foundation. G.S.\ thanks the Deutsche
Forschungsgemeinschaft and the Gesellschaft f\"{u}r Schwerionenforschung
for support. E.J.W.\ thanks Prof.\ N.\ Temme for providing information
on the Van Wijngaarden transformation and Prof.\ A.\ Z.\ Msezane and
Prof.\ C.\ R.\ Handy for their invitation to the Center for Theoretical
Studies of Physical Systems at Clark Atlanta University. E.J.W.\
acknowledges support from the Fonds der Chemischen Industrie.

\appendix

\section{On the Efficiency of the Transformation of Alternating
and Monotone Series}
\label{App_A}

In the case of the Levin transformation ${\cal L}_{k}^{(n)} (\beta, s_n,
\omega_n)$, Eq.\ (\ref{LevTr}), and the closely related Weniger
transformation ${\cal S}_{k}^{(n)} (\beta, s_n, \omega_n)$, Eq.
(\ref{SidTr}), it is relatively easy to understand why alternating
series can be transformed much more effectively than monotone series.
These two transformations are both special cases of the following
transformation which is characterized by the positive weights $w_k (n)$:
\begin{equation}
{\cal T}_k^{(n)} \bigl( w_k (n); s_n, \omega_n \bigr)
\; = \; \frac
{\Delta^k \{ w_k (n) s_n / \omega_n \} }
{\Delta^k \{ w_k (n) / \omega_n \} } \, .
\label{GenTr}
\end{equation}
If $w_k (n) = (n + \beta)^{k-1}$, we obtain Levin's transformation, and
it $w_k (n) = (n + \beta)_{k-1}$, we obtain Weniger's transformation.

Since the difference operator $\Delta^k$ is linear, Eq.\ (\ref{GenTr})
can also be rewritten as follows:
\begin{eqnarray}
{\cal T}_k^{(n)} \bigl( w_k (n); s_n, \omega_n \bigr)
& = & s \, + \, \frac
{\Delta^k \{ w_k (n) [s_n-s] / \omega_n \} }
{\Delta^k \{ w_k (n) / \omega_n \} } \, .
\label{T_Rat}
\end{eqnarray}
Obviously, the sequence transformation ${\cal T}_k^{(n)}$ converges to
the (generalized) limit $s$ of the input sequence $\{ s_n
\}_{n=0}^{\infty}$ if the remainder estimates $\{ \omega_n
\}_{n=0}^{\infty}$ can be chosen in such a way that the ratio on the
right-hand side becomes negligibly small.

Loosely speaking, this means that we have to choose the remainder
estimates in such a way that the numerator of the ratio on the
right-hand side of Eq.\ (\ref{T_Rat}) becomes as small as possible and
the denominator becomes as large as possible.

If we can find remainder estimates such that
\begin{equation}
s_n - s \; = \; \omega_n \bigl[ c + O (n^{-1}) \bigr] \, ,
\qquad n \to \infty \, ,
\label{asy_cond}
\end{equation}
then the weighted difference operator $\Delta^k w_k (n)$ will make the
numerator small, no matter whether the input data $\{
s_n\}_{n=0}^{\infty}$ are the partial sums of an alternating or of a
monotone series.

The situation is quite different in the case of the denominator.
Application of the weighted difference operator $\Delta^k w_k (n)$ to
$1/\omega_n$ yields
\begin{eqnarray}
\lefteqn{\Delta^k \{ w_k (n) / \omega_n \}} \nonumber \\
& \qquad = & (-1)^k \,
\sum_{j=0}^{k} \, (-1)^j \, {{k} \choose {j}} \,
\frac {w_k (n+j)} {\omega_{n+j}} \, .
\end{eqnarray}
If the remainder estimates all have the same sign, the denominator will
also become relatively small because of cancellation due to differencing
(hopefully not as small as the numerator because then the transformation
process would not converge). This cancellation process is also the major
source of numerical instabilities which are magnified since they occur
in the denominator.

If the remainder estimates are strictly alternating in sign, $\omega_n =
(-1)^n \vert \omega_n \vert$, then we obtain ($w_k (n)$ is by assumption
positive)
\begin{equation}
\Delta^k \{ w_k (n)/\omega_n \} \; = \; (-1)^{k+n} \,
\sum_{j}^{k} \, {{k} \choose {j}} \,
\frac {w_k (n+j)} {\vert \omega_{n+j} \vert} \, .
\label{DenomSum}
\end{equation}
Thus, all terms of the denominator sum have the same sign, and there is
no cancellation as in the case of monotone remainder estimates and also
no source of numerical instabilities. Consequently, the denominator
becomes relatively large which improves convergence.

As an example, we apply the transformation (\ref{GenTr}) to the monotone
model sequence
\begin{equation}
s_n \; = \; s \, + \, \frac {c_0} {(n+1)^{\alpha}} \, + \,
\frac {c_1} {(n+1)^{\alpha+1}} \, + \, \ldots \, ,
\label{MonModSeq}
\end{equation}
as well as to the alternating model sequence
\begin{eqnarray}
\lefteqn{t_n \; = \; t} \nonumber \\
& & \qquad + \, (-1)^n \, \left\{ \frac {c_0} {(n+1)^{\alpha}}
\, + \, \frac {c_1} {(n+1)^{\alpha+1}} \, + \, \ldots \right\}
\label{AltModSeq}
\end{eqnarray}
and derive asymptotic ($n \to \infty$) transformation error estimates.

For that purpose, we assume $\alpha > 0$ -- which implies that the two
model sequences converge to their limits $s$ and $t$, respectively --
and $w_k (n) = O \bigl( n^{k-1} \bigr)$ as $n \to \infty$. In the case
of the monotone input sequence $\{ s_n\}_{n=0}^{\infty}$, we choose
$\omega_n = (n+1)^{-\alpha}$, and in the case of the alternating input
sequence $\{ t_n \}_{n=0}^{\infty}$, we choose $\omega_n = (-1)^n
(n+1)^{-\alpha}$. Then, we obtain for the numerators
\begin{eqnarray}
\lefteqn{\Delta^k \left\{
\frac {w_k (n) [s_n-s]} {(n+1)^{-\alpha}} \right\}}
\nonumber \\
& \qquad = & \Delta^k \left\{ \frac {w_k (n) [t_n-t]}
{(-1)^n (n+1)^{-\alpha}} \right\}
\; = \; O \bigl( n^{-k-1} \bigr) \, .
\label{NumEst}
\end{eqnarray}
The estimates for the denominators differ considerably. In the case of
the monotone input sequence (\ref{MonModSeq}), we obtain
\begin{equation}
\Delta^k \left\{
\frac {w_k (n)} {(n+1)^{-\alpha}} \right\} \; = \;
O \bigl( n^{\alpha-1} \bigr) \, ,
\end{equation}
which in combination with Eq.\ (\ref{NumEst}) yields the following
asymptotic transformation error estimate as $n \to \infty$:
\begin{equation}
\frac
{{\cal T}_k^{(n)} \bigl( w_k (n); s_n, (n+1)^{-\alpha} \bigr) - s}
{s_n - s} \; = \; O \bigl( n^{-k} \bigr) \, .
\label{EstMon}
\end{equation}
In the case of the alternating input sequence (\ref{AltModSeq}), we
exploit the fact that the first term with $j = 0$ in the sum on the
right-hand side of Eq.\ (\ref{DenomSum}) is smaller in magnitude than
the whole sum, yielding
\begin{eqnarray}
\lefteqn{\frac {1} {\left\vert
\Delta^k \{ w_k (n)/[(-1)^n (n+1)^{-\alpha}] \} \right\vert}}
\nonumber \\
& \qquad \le &
\frac {1}{\left\vert w_k (n)/[(-1)^n (n+1)^{-\alpha}] \right\vert} \, .
\label{DenomEst}
\end{eqnarray}
Of course, we could also try to construct more sophisticated estimates
for the denominators. However, the relatively crude estimate
(\ref{DenomEst}) suffices for our purposes since it implies in
combination with Eq.\ (\ref{NumEst})
\begin{equation}
\frac
{{\cal T}_k^{(n)} \bigl( w_k (n); t_n,
(-1)^n (n+1)^{-\alpha} \bigr) - t}
{t_n - t} \; = \; O \bigl( n^{-2 k} \bigr) \, .
\label{EstAlt}
\end{equation}
The two transformation error estimates (\ref{EstMon}) and
(\ref{EstAlt}), which both hold as $n \to \infty$, show that there is a
substantial difference between the transformation of monotone and
alternating series. There is a considerable amount of evidence that this
conclusion is actually generally true, i.e., also in the case of Pad\'e
approximants and other sequence transformations.

\section{Exactness results}
\label{App_B}

In this Appendix, we analyze exactness properties of the Levin
transformation ${\cal L}_{k}^{(n)} (\beta, s_n, \omega_n)$, Eq.\
(\ref{LevTr}), and of the Weniger transformation ${\cal S}_{k}^{(n)}
(\beta, s_n, \omega_n)$, Eq.\ (\ref{SidTr}), and their variants
$d_{k}^{(n)} \bigl( \beta, {\bf S}_n \bigr)$, Eq.\ (\ref{dLevTr}), and
${\delta}_{k}^{(n)} \bigl( \beta, {\bf S}_n \bigr)$, Eq.\
(\ref{dSidTr}), which both use the remainder estimate (\ref{d_Est})
proposed by Smith and Ford \cite{SmFo79}.

For that purpose we introduce in Eq.\ (\ref{T_Rat}) the remainder $r_n =
s_n - s$ and obtain
\begin{eqnarray}
{\cal T}_k^{(n)} \bigl( w_k (n); s_n, \omega_n \bigr)
& = & s \, + \, \frac
{\Delta^k \{ w_k (n) r_n / \omega_n \} }
{\Delta^k \{ w_k (n) / \omega_n \} } \, .
\label{T_Rat_rem}
\end{eqnarray}
Obviously, the general sequence transformation ${\cal T}_k^{(n)}$, Eq.\
(\ref{GenTr}), which contains ${\cal L}_{k}^{(n)} (\beta, s_n,
\omega_n)$ and ${\cal S}_{k}^{(n)} (\beta, s_n, \omega_n)$ as special
cases, is exact for a given input sequence $\{ s_n \}_{n=0}^{\infty}$ if
the difference operator $\Delta^k$ annihilates $w_k (n) r_n / \omega_n$
but not $w_k (n) / \omega_n$.

Thus, if we want to prove the exactness of ${\cal T}_k^{(n)}$ for some
sequence $\{ s_n \}_{n=0}^{\infty}$, we need to know explicit
expressions for the remainders $\{ r_n \}_{n=0}^{\infty}$ of the input
sequence as functions of $n$.

In the case of the alternating series (\ref{ZetaAlt}) for the Riemann
zeta function with zero or negative integral argument, this can be
accomplished easily. Replacing $z$ by $- l$ in Eq.\ (\ref{ZetaAlt}) with
$l = 0, 1, 2, \ldots$ we obtain
\begin{equation}
\zeta (-l) \; = \; \frac {1}{1-2^{l+1}} \,
\sum_{j=0}^{\infty} \, (-1)^j (j+1)^l \, .
\label{ZetaNegArg}
\end{equation}
For $l = 0$, this series is the negative of the geometric series
$1/(1+x) = \sum_{j=0}^{\infty} (-x)^j$ at $x = 1$. According to Theorem
12-9 of Ref.\ \cite{We89}, the sequence transformations ${\cal
L}_{k}^{(n)} (\beta, s_n, \omega_n)$ and ${\cal S}_{k}^{(n)} (\beta,
s_n, \omega_n)$ are exact for the geo\-metric series if the first term
neglected in the partial sum is used as remainder estimate (which
corresponds to the remainder estimate (\ref{d_Est})). Thus, the
exactness of these sequence transformations for the infinite series
(\ref{ZetaNegArg}) has to be analyzed only for $l \ge 1$.

For the determination of its truncation error with arbitrary integral $l
\ge 1$, we rewrite the infinite series (\ref{ZetaNegArg}) as follows:
\begin{eqnarray}
\lefteqn{\sum_{j=0}^{\infty} \, (-1)^j (j+1)^l} \nonumber \\
& & \quad = \; \sum_{j=0}^{n} \, (-1)^j (j+1)^l \, + \,
\sum_{j=n+1}^{\infty} \, (-1)^j (j+1)^l
\\
& & \quad = \; \sum_{j=0}^{n} \, (-1)^j (j+1)^l \nonumber \\
& & \qquad + \,
(-1)^{n+1} \, \sum_{\nu=0}^{\infty} \, (-1)^{\nu} (n+\nu+2)^l \, .
\label{TruncErrRep}
\end{eqnarray}
The infinite series on the right-hand side of Eq.\ (\ref{TruncErrRep})
obviously diverges. However, it can be summed easily since it can be
represented as a derivative of the geometric series $1/(1+x) =
\sum_{j=0}^{\infty} (-x)^j$ at $x = 1$:
\begin{eqnarray}
\lefteqn{\sum_{\nu=0}^{\infty} \, (-1)^{\nu} (n+\nu+2)^l} \nonumber \\
& & \quad = \; \lim_{x \to 1-} \, \left\{
\sum_{\nu=0}^{\infty} \, (-1)^{\nu} (n+\nu+2)^l x^{n+\nu+1} \right\}
\\
& & \quad = \; \lim_{x \to 1-} \, \left\{
\left( \frac {{\rm d}}{{\rm d} x} x \right)^l \,
\sum_{\nu=0}^{\infty} \, (-1)^{\nu} x^{n+\nu+1} \right\}
\\
& & \quad = \; \left\{
\left( \frac {{\rm d}}{{\rm d} x} x \right)^l \,
\frac {x^{n+1}}{1+x} \right\} \Bigg\vert_{x=1} \, .
\label{ZetaNegArgTruncErr}
\end{eqnarray}
By inserting this into Eq.\ (\ref{ZetaNegArg}) we obtain the following
representation of $\zeta (-l)$ as a partial sum plus an explicit
expression for the truncation error:
\begin{eqnarray}
\zeta (-l) & & \quad = \; \frac {1}{1-2^{l+1}} \,
\Biggl\{ \sum_{j=0}^{n} \, (-1)^j (j+1)^l
\nonumber \\
& & \qquad + \, (-1)^{n+1} \,
\left[ \left( \frac {{\rm d}}{{\rm d} x} x \right)^l \,
\frac {x^{n+1}}{1+x} \right] \Bigg\vert_{x=1} \Biggr\} \, .
\label{ZetaNegArgPSplusTruncErr}
\end{eqnarray}
If we set in Eq.\ (\ref{ZetaNegArgPSplusTruncErr}) $n = - 1$, then the
partial sum is an empty sum and vanishes, yielding the following
explicit expression for the Riemann zeta function with zero or negative
integer argument:
\begin{equation}
\zeta (-l) \; = \; \frac {1}{1-2^{l+1}} \, \left[
\left( \frac {{\rm d}}{{\rm d} x} x \right)^l \,
\frac {1}{1+x} \right] \Bigg\vert_{x=1} \, .
\label{ZetaNegArgExplExpr}
\end{equation}
Equivalent explicit expressions can be obtained via the substitution $y
= \ln (x)$. Obviously, we have
\begin{equation}
\left[ \frac {{\rm d}} {{\rm d} y} + 1 \right] f \bigl({\rm e}^y \bigr)
\bigg\vert_{y=\ln(x)} \; = \;
\frac {{\rm d}} {{\rm d} x} \bigr[ x \, f (x) \bigr] \, .
\end{equation}
In this way, we obtain:
\begin{eqnarray}
\lefteqn{\zeta (-l)} \nonumber \\
& & \quad = \; \frac {1}{1-2^{l+1}} \, \left\{
\left[ \frac {{\rm d}}{{\rm d} y} + 1 \right]^l \,
\frac {1}{1+{\rm e}^y} \right\} \Bigg\vert_{y=0}
\label{ExpoZetaNegArgExplExpr}
\\
& & \quad = \; \frac {1}{1-2^{l+1}} \, \Biggl\{ \sum_{j=0}^{n} \,
(-1)^j (j+1)^l
\nonumber \\
& & \qquad \; + \, (-1)^{n+1} \,
\left[ \left( \frac {{\rm d}}{{\rm d} y} + 1 \right)^l \,
\frac {{\rm e}^{(n+1)y}}{1+{\rm e}^y} \right] \Bigg\vert_{y=0}
\Biggr\} \, .
\end{eqnarray}
Combination of Eqs.\ (\ref{BernoulZetaNegArg}),
(\ref{ZetaNegArgExplExpr}), and (\ref{ExpoZetaNegArgExplExpr}) yields
the following expressions for the Bernoulli numbers with even indices:
\begin{eqnarray}
B_{2 l} & = & \frac {- 2 l}{1-2^{2 l}} \, \left[
\left( \frac {{\rm d}}{{\rm d} x} x \right)^{2l-1} \,
\frac {1}{1+x} \right] \Bigg\vert_{x=1}
\\
& = & \frac {- 2 l}{1-2^{2 l}} \, \left\{
\left[ \frac {{\rm d}}{{\rm d} y} + 1 \right]^{2l-1} \,
\frac {1}{1+{\rm e}^y} \right\} \Bigg\vert_{y=0} \, .
\end{eqnarray}

It follows from Eq.\ (\ref{T_Rat_rem}) that the general sequence
transformation ${\cal T}_k^{(n)}$, Eq.\ (\ref{GenTr}), is exact for some
input sequence $\{ s_n \}_{n=0}^{\infty}$ if the expression $w_k (n) r_n
/ \omega_n$ with $r_n = s_n - s$ is a polynomial of degree $k-1$ in $n$.
Thus, Eq.\ (\ref{ZetaNegArgTruncErr}) implies that in the case of the
Weniger transformation ${\delta}_{n}^{(0)} \bigl( 1, {\bf S}_0 \bigr)$,
Eq.\ (\ref{dSidTr}), which uses the remainder estimate (\ref{d_Est}), we
have to show that the ratio
\begin{equation}
\frac {(n+1)_{k-1} \sum_{\nu=0}^{\infty} \, (-1)^{\nu} (n+\nu+2)^l}
{(-1)^{n+1} \, (n+2)^l}
\label{WeRatZeta}
\end{equation}
is for sufficiently large values of $k$ a polynomial of degree $k-1$ in
$n$. For $l = 1$, 2, 3, and 4, Eq.\ (\ref{ZetaNegArgTruncErr}) yields
the following explicit expressions:
\begin{eqnarray}
\lefteqn{\sum_{\nu=0}^{\infty} \, (-1)^{\nu} (n+\nu+2)^1 \; = \;
\frac{2 n + 3}{4} \, ,} 
\\
\lefteqn{\sum_{\nu=0}^{\infty} \, (-1)^{\nu} (n+\nu+2)^2 \; = \;
\frac{(n + 1) (n + 2)}{2} \, ,} 
\\
\lefteqn{\sum_{\nu=0}^{\infty} \, (-1)^{\nu} (n+\nu+2)^3}
\nonumber \\
& & \qquad = \; \frac{(2 n + 3) (2 n^2 + 6 n + 3)}{8} \, ,
\\
\lefteqn{\sum_{\nu=0}^{\infty} \, (-1)^{\nu} (n+\nu+2)^4}
\nonumber \\
& & \qquad = \; \frac{(n + 1) (n + 2) (n^2 + 3 n + 1)}{2} \, .
\end{eqnarray}
Since $(n+1)_{k-1} = (n+1)(n+2) \ldots (n+k-1)$, the ratio
(\ref{WeRatZeta}) is for $l = 1$, 2 and $k \ge 3$ a polynomial of degree
$k-1$ in $n$, whereas for $l = 3$, 4 it is a rational expression in $n$
which is not annihilated by the operator $\Delta^k$.

Thus, the exactness of the transformation ${\delta}_{n}^{(0)} \bigl(1,
{\bf S}_0 \bigr)$ for $\zeta (-1)$ as in Table \ref{Tab_4_2}, or for
$\zeta (-2)$ is more or less accidental. If, however, we choose in the
Levin transformation $d_{n}^{(0)} \bigl(\beta, {\bf S}_0 \bigr)$, Eq.\
(\ref{dLevTr}), with $\beta = 2$ instead of the usual $\beta = 1$, we
find that the corresponding ratio
\begin{equation}
\frac {(n+2)^{k-1} \sum_{\nu=0}^{\infty} \, (-1)^{\nu} (n+\nu+2)^l}
{(-1)^{n+1} \, (n+2)^l}
\label{LevRatZeta}
\end{equation}
is a polynomial of degree $k - 1$ in $n$ for all $l = 0$, $1$, $2$,
$\ldots$ and $k \ge l + 1$. This follows at once from the fact that the
differential operator in Eq.\ (\ref{ZetaNegArgTruncErr}) produces a
polynomial of degree $l$ in $n$.

In the same way, it can be shown that the Weniger transformation
${\delta}_{n}^{(0)} \bigl( \beta, {\bf S}_0 \bigr)$, Eq.\
(\ref{dSidTr}), with $\beta = 2$ is exact for the hypergeometric series
for all $l = 0$, $1$, $2$, $\ldots$ and for $k \ge l + 1$:
\begin{equation}
l! {}_1 F_0 (l+1; -1) \; = \;
\sum_{j=0}^{\infty} \, (-1)^j (j+1)_l \; = \;
\frac {l!} {2^{l+1}} \, .
\end{equation}
This can also be rewritten as follows:
\begin{eqnarray}
\lefteqn{l! {}_1 F_0 (l+1; -1)} \nonumber \\
& & = \; \sum_{j=0}^{n} \, (-1)^j (j+1)_l \nonumber \\
& & \quad + \,
(-1)^{n+1} \, \sum_{\nu=0}^{\infty} \, (-1)^{\nu} (n+\nu+2)_l
\\
& & = \; \sum_{j=0}^{n} \, (-1)^j (j+1)_l \nonumber \\
& & \quad + \, (-1)^{n+1} \, (n+2)_l \, {}_2 F_1 (1, n+l+2; n+2; -1)
\\
& & = \; \sum_{j=0}^{n} \, (-1)^j (j+1)_l \nonumber \\
& & \quad + \, (-1)^{n+1} \,
\sum_{j=0}^{l} \, (-l)_j \, (n+j+2)_{l-j} \, (1/2)^{l+1} \, .
\label{1F0_RemEst}
\end{eqnarray}
In this case, the corresponding ratio
\begin{equation}
\frac {(n+2)_{k-1} \sum_{\nu=0}^{\infty} \, (-1)^{\nu} (n+\nu+2)_l}
{(-1)^{n+1} \, (n+2)_l}
\end{equation}
is a polynomial of degree $k - 1$ in $n$ for all $l = 0$, 1, 2, $\ldots$
and for $k \ge l + 1$. This follows at once from the fact that the
second sum on the right-hand side of Eq.\ (\ref{1F0_RemEst}) is a
polynomial of degree $l$ in $n$.

However, the value $\beta = 2$ in the Eqs.\ (\ref{dLevTr}) and
(\ref{dSidTr}) cannot significantly improve convergence for the zeta
function. We also investigated the performance of the transformations
$d_{n}^{(0)} \bigl( 2, {\bf S}_n)$ and ${\delta}_{n}^{(0)} \bigl( 2,
{\bf S}_n)$. Except for the cases of $z = -1$ or $z = -2$ in which the
Levin transformation becomes exact with $\beta = 2$, the transformations
$d_{n}^{(0)} \bigl( 2, {\bf S}_n)$ and ${\delta}_{n}^{(0)} \bigl( 2,
{\bf S}_n)$ yield results which differ only marginally from the results
obtained by $d_{n}^{(0)} \bigl( 1, {\bf S}_n)$ and ${\delta}_{n}^{(0)}
\bigl( 1, {\bf S}_n)$, which are presented in Tables \ref{Tab_4_1} -
\ref{Tab_4_3}.

\newpage

\onecolumn
\centerline{\bf \large TABLES}

\mediumtext

\begin{table}
\caption{Demonstration that Van Wijngaarden's transformation, Eqs.\
(\protect\ref{VWSerTran}) - (\protect\ref{B2a}), is a rearrangement of
the original series $\sum_{k=0}^{\infty} a(k)$.}
\label{Tab_2_1}
\begin{tabular}{c|rrrrrlr}
${(-1)^j {\bf A}}_j$ & \multicolumn{7}{c}
{Leading terms of the original series} \\
\tableline \rule[-4pt]{0pt}{3.8\jot}%
$\phantom{-}
  {\bf A}_0$ & $a(0)$ & $+ 2 \, a (1)$ & & $+ 4 \, a(3)$ & & $+$ &
  $\dots$ \\
$-{\bf A}_1$ & & $- a(1)$ & & $- 2 \, a(3)$ & & $-$ & $\dots$ \\
$\phantom{-}
  {\bf A}_2$  & & & $+a(2)$ & & & $+$ & $\dots$ \\
$-{\bf A}_3$ & & & & $- a(3)$ & & $-$ & $\dots$ \\
$\phantom{-}
  {\bf A}_4$  & & & & & $+ a(4)$ & $-$ & $\dots$ \\
\tableline \rule{0pt}{3.5\jot}%
$\sum_{j=0}^{\infty} (-1)^j {\bf A}_j$ & $a(0)$ & $+ a(1)$ & $+ a(2)$
& $+ a(3)$ & $+ a(4)$ & $+$ & $\dots$ \\
\end{tabular}
\end{table}

\widetext

\begin{table}
\caption{Evaluation of the Riemann zeta function of argument $z = 1.01$
by accelerating the convergence of the alternating series
(\protect\ref{ZetaAlt}). The result is given as $10^{-3}\,\zeta(1.01)$.}
\label{Tab_4_1}
\begin{tabular}{lrrrr}%
$n$%
& \multicolumn{1}{c}{${\bf S}_{n}$}%
& \multicolumn{1}{c}{${\bf E}_n$}%
& \multicolumn{1}{c}{$d_{n}^{(0)} \bigl(1, {\bf S}_0 \bigr)$}%
& \multicolumn{1}{c}{${\delta}_{n}^{(0)} \bigl(1, {\bf S}_0 \bigr)$}%
\\
\tableline \rule[-4pt]{0pt}{3.8\jot}%
 0 & 0.144~770~081~711~084 & 0.144~770~081~711~084 & 0.144~770~081~711~084 & 0.144~770~081~711~084 \\
 1 & 0.072~885~040~855~542 & 0.090~606~301~069~428 & 0.101~569~133~143~252 & 0.101~569~133~143~252 \\
 2 & 0.120~614~482~322~980 & 0.096~697~481~252~857 & 0.100~456~642~533~059 & 0.100~456~642~533~059 \\
 3 & 0.084~920~235~019~068 & 0.098~985~546~018~036 & 0.100~587~783~459~042 & 0.100~579~332~613~649 \\
 4 & 0.113~411~984~373~518 & 0.099~901~837~494~047 & 0.100~577~428~866~203 & 0.100~578~083~572~921 \\
 5 & 0.089~712~109~307~161 & 0.100~283~957~662~399 & 0.100~577~954~415~585 & 0.100~577~949~566~834 \\
 6 & 0.109~994~997~614~328 & 0.100~447~835~715~031 & 0.100~577~944~116~204 & 0.100~577~943~567~122 \\
 7 & 0.092~271~153~050~434 & 0.100~519~572~572~454 & 0.100~577~943~249~050 & 0.100~577~943~346~150 \\
 8 & 0.108~007~136~313~467 & 0.100~551~470~653~282 & 0.100~577~943~342~049 & 0.100~577~943~338~734 \\
 9 & 0.093~859~665~080~579 & 0.100~565~830~599~811 & 0.100~577~943~338~553 & 0.100~577~943~338~503 \\
10 & 0.106~708~750~240~925 & 0.100~572~360~153~981 & 0.100~577~943~338~482 & 0.100~577~943~338~497 \\
11 & 0.094~940~666~205~327 & 0.100~575~353~801~276 & 0.100~577~943~338~498 & 0.100~577~943~338~497 \\
12 & 0.105~794~821~569~601 & 0.100~576~735~870~616 & 0.100~577~943~338~497 & 0.100~577~943~338~497 \\
13 & 0.095~723~429~487~784 & 0.100~577~377~707~430 & 0.100~577~943~338~497 & 0.100~577~943~338~497 \\
14 & 0.105~116~912~361~076 & 0.100~577~677~299~954 & 0.100~577~943~338~497 & 0.100~577~943~338~497 \\
15 & 0.096~316~203~847~728 & 0.100~577~817~763~434 & 0.100~577~943~338~497 & 0.100~577~943~338~497 \\
\tableline \rule{0pt}{3.5\jot}%
exact & 0.100~577~943~338~497 & 0.100~577~943~338~497 & 0.100~577~943~338~497 & 0.100~577~943~338~497
\end{tabular}
\end{table}

\mediumtext

\begin{table}
\caption{Evaluation of the Riemann zeta function of argument $z = - 1$
by summing the divergent alternating series (\protect\ref{ZetaAlt}).
The result is given as $10\,\zeta(-1)$.}
\label{Tab_4_2}
\begin{tabular}{lrrr}
$n$ &
\multicolumn{1}{c}{${\bf S}_{n}$} &
\multicolumn{1}{c}{$d_{n}^{(0)} \bigl(1, {\bf S}_0 \bigr)$} &
\multicolumn{1}{c}{${\delta}_{n}^{(0)} \bigl(1, {\bf S}_0 \bigr)$} \\
\tableline \rule[-4pt]{0pt}{3.8\jot}%
 0 & $ -3.333~333~333~333~333$ & $-3.333~333~333~333~333$ & $-3.333~333~333~333~333$ \\
 1 & $  3.333~333~333~333~333$ & $-0.666~666~666~666~667$ & $-0.666~666~666~666~667$ \\
 2 & $ -6.666~666~666~666~667$ & $-0.860~215~053~763~441$ & $-0.860~215~053~763~441$ \\
 3 & $  6.666~666~666~666~667$ & $-0.830~449~826~989~619$ & $-0.833~333~333~333~333$ \\
 4 & $-10.000~000~000~000~000$ & $-0.833~557~890~954~819$ & $-0.833~333~333~333~333$ \\
 5 & $ 10.000~000~000~000~000$ & $-0.833~319~627~418~409$ & $-0.833~333~333~333~333$ \\
 6 & $-13.333~333~333~333~333$ & $-0.833~334~020~666~741$ & $-0.833~333~333~333~333$ \\
 7 & $ 13.333~333~333~333~333$ & $-0.833~333~304~093~535$ & $-0.833~333~333~333~333$ \\
 8 & $-16.666~666~666~666~667$ & $-0.833~333~334~413~139$ & $-0.833~333~333~333~333$ \\
 9 & $ 16.666~666~666~666~667$ & $-0.833~333~333~298~109$ & $-0.833~333~333~333~333$ \\
10 & $-20.000~000~000~000~000$ & $-0.833~333~333~334~362$ & $-0.833~333~333~333~333$ \\
11 & $ 20.000~000~000~000~000$ & $-0.833~333~333~333~306$ & $-0.833~333~333~333~333$ \\
12 & $-23.333~333~333~333~333$ & $-0.833~333~333~333~334$ & $-0.833~333~333~333~333$ \\
13 & $ 23.333~333~333~333~330$ & $-0.833~333~333~333~333$ & $-0.833~333~333~333~333$ \\
14 & $-26.666~666~666~666~670$ & $-0.833~333~333~333~333$ & $-0.833~333~333~333~333$ \\
15 & $ 26.666~666~666~666~670$ & $-0.833~333~333~333~333$ & $-0.833~333~333~333~333$ \\
\hline
\tableline \rule{0pt}{3.5\jot}%
exact & $-0.833~333~333~333~333$ & $-0.833~333~333~333~333$ & $-0.833~333~333~333~333$
\end{tabular}
\end{table}

\widetext

\begin{table}
\caption{Evaluation of $\zeta(1/2 + 13.7~i)$ with the CNCT.}%
\label{Tab_4_3}%
\begin{tabular}{lrcrrcr}%
$n$%
& \multicolumn{3}{c}{$d_{n}^{(0)} \bigl(1, {\bf S}_0\bigr)$}%
& \multicolumn{3}{c}{${\delta}_{n}^{(0)} \bigl(1, {\bf S}_0\bigr)$}%
\\
\tableline \rule[-4pt]{0pt}{3.8\jot}%
    0 & $0.414~107~543~949~134$ & $+$ & $0.017~316~297~125~790~i$ & $ 0.414~107~543~949~134$ & $+$ & $0.017~316~297~125~790~i$ \\
    1 & $0.575~871~239~097~112$ & $-$ & $0.042~690~435~565~758~i$ & $ 0.575~871~239~097~112$ & $-$ & $0.042~690~435~565~758~i$ \\
    2 & $0.523~424~912~174~020$ & $+$ & $0.152~835~043~959~961~i$ & $ 0.567~958~887~269~553$ & $+$ & $0.129~913~386~486~220~i$ \\
    3 & $0.481~953~715~196~159$ & $-$ & $0.288~400~086~031~923~i$ & $ 0.474~129~917~411~618$ & $-$ & $0.168~775~423~138~291~i$ \\
    4 & $0.012~442~899~246~184$ & $-$ & $0.237~603~260~694~125~i$ & $-0.180~827~868~994~142$ & $-$ & $0.367~542~940~737~051~i$ \\
    5 & $0.123~074~021~609~358$ & $-$ & $0.316~357~718~264~423~i$ & $ 0.126~392~529~409~594$ & $-$ & $0.290~235~127~404~228~i$ \\
    6 & $0.105~377~569~236~175$ & $-$ & $0.313~246~538~788~829~i$ & $ 0.107~386~234~787~298$ & $-$ & $0.317~087~400~005~856~i$ \\
    7 & $0.107~635~288~132~001$ & $-$ & $0.312~843~180~304~712~i$ & $ 0.107~124~241~668~490$ & $-$ & $0.312~622~662~817~549~i$ \\
    8 & $0.107~429~326~957~578$ & $-$ & $0.312~999~708~812~577~i$ & $ 0.107~489~121~873~498$ & $-$ & $0.312~978~336~617~038~i$ \\
    9 & $0.107~438~933~679~469$ & $-$ & $0.312~974~184~308~675~i$ & $ 0.107~436~183~812~222$ & $-$ & $0.312~980~180~873~835~i$ \\
   10 & $0.107~439~640~888~613$ & $-$ & $0.312~976~813~188~762~i$ & $ 0.107~439~393~557~558$ & $-$ & $0.312~976~229~866~877~i$ \\
   11 & $0.107~439~434~156~190$ & $-$ & $0.312~976~661~398~796~i$ & $ 0.107~439~488~450~447$ & $-$ & $0.312~976~678~675~422~i$ \\
   12 & $0.107~439~457~152~512$ & $-$ & $0.312~976~659~181~027~i$ & $ 0.107~439~453~065~793$ & $-$ & $0.312~976~661~841~904~i$ \\
   13 & $0.107~439~455~840~151$ & $-$ & $0.312~976~660~717~255~i$ & $ 0.107~439~455~877~677$ & $-$ & $0.312~976~660~319~609~i$ \\
   14 & $0.107~439~455~825~365$ & $-$ & $0.312~976~660~547~552~i$ & $ 0.107~439~455~848~368$ & $-$ & $0.312~976~660~568~977~i$ \\
   15 & $0.107~439~455~836~355$ & $-$ & $0.312~976~660~556~014~i$ & $ 0.107~439~455~833~989$ & $-$ & $0.312~976~660~556~440~i$ \\
   16 & $0.107~439~455~835~269$ & $-$ & $0.312~976~660~556~233~i$ & $ 0.107~439~455~835~348$ & $-$ & $0.312~976~660~556~072~i$ \\
   17 & $0.107~439~455~835~311$ & $-$ & $0.312~976~660~556~157~i$ & $ 0.107~439~455~835~317$ & $-$ & $0.312~976~660~556~169~i$ \\
   18 & $0.107~439~455~835~313$ & $-$ & $0.312~976~660~556~163~i$ & $ 0.107~439~455~835~312$ & $-$ & $0.312~976~660~556~163~i$ \\
   19 & $0.107~439~455~835~313$ & $-$ & $0.312~976~660~556~163~i$ & $ 0.107~439~455~835~313$ & $-$ & $0.312~976~660~556~163~i$ \\
   20 & $0.107~439~455~835~313$ & $-$ & $0.312~976~660~556~163~i$ & $ 0.107~439~455~835~313$ & $-$ & $0.312~976~660~556~163~i$ \\
   21 & $0.107~439~455~835~313$ & $-$ & $0.312~976~660~556~163~i$ & $ 0.107~439~455~835~313$ & $-$ & $0.312~976~660~556~163~i$ \\
   22 & $0.107~439~455~835~313$ & $-$ & $0.312~976~660~556~163~i$ & $ 0.107~439~455~835~313$ & $-$ & $0.312~976~660~556~163~i$ \\
   23 & $0.107~439~455~835~313$ & $-$ & $0.312~976~660~556~163~i$ & $ 0.107~439~455~835~313$ & $-$ & $0.312~976~660~556~163~i$ \\
   24 & $0.107~439~455~835~313$ & $-$ & $0.312~976~660~556~163~i$ & $ 0.107~439~455~835~313$ & $-$ & $0.312~976~660~556~163~i$ \\
   25 & $0.107~439~455~835~313$ & $-$ & $0.312~976~660~556~163~i$ & $ 0.107~439~455~835~313$ & $-$ & $0.312~976~660~556~163~i$ \\
\tableline \rule{0pt}{3.5\jot}%
exact & $0.107~439~455~835~313$ & $-$ & $0.312~976~660~556~163~i$ & $ 0.107~439~455~835~313$ & $-$ & $0.312~976~660~556~163~i$
\end{tabular}
\end{table}

\mediumtext

\begin{table}
\caption{Evaluation of
$10^{-2}\,{\rm Li}_1 (0.99999) = - 10^{-2}\,\ln (0.00001)$ with the CNCT.}
\label{Tab_5_1}%
\begin{tabular}{lrrr}%
$n$%
& \multicolumn{1}{c}{${\bf S}_{n}$}%
& \multicolumn{1}{c}{$d_{n}^{(0)} \bigl(1, {\bf S}_0 \bigr)$}%
& \multicolumn{1}{c}{${\delta}_{n}^{(0)} \bigl(1, {\bf S}_0 \bigr)$}%
\\
\tableline \rule[-4pt]{0pt}{3.8\jot}%
 0 & 0.162~768~973~713~089 & 0.162~768~973~713~089 & 0.162~768~973~713~089 \\
 1 & 0.086~384~436~856~544 & 0.116~225~388~785~336 & 0.116~225~388~785~336 \\
 2 & 0.135~357~615~351~010 & 0.114~995~657~006~664 & 0.114~995~657~006~664 \\
 3 & 0.099~665~296~922~988 & 0.115~140~148~148~939 & 0.115~131~002~772~470 \\
 4 & 0.127~575~317~431~702 & 0.115~128~665~188~679 & 0.115~129~400~970~919 \\
 5 & 0.104~755~344~851~635 & 0.115~129~271~216~942 & 0.115~129~261~517~373 \\
 6 & 0.123~997~631~730~577 & 0.115~129~254~851~082 & 0.115~129~254~830~226 \\
 7 & 0.107~401~422~517~315 & 0.115~129~254~746~908 & 0.115~129~254~725~403 \\
 8 & 0.121~964~815~378~398 & 0.115~129~254~612~355 & 0.115~129~254~664~668 \\
 9 & 0.109~009~755~125~042 & 0.115~129~254~647~759 & 0.115~129~254~647~723 \\
10 & 0.120~662~087~466~136 & 0.115~129~254~650~661 & 0.115~129~254~648~068 \\
11 & 0.110~085~384~510~686 & 0.115~129~254~650~507 & 0.115~129~254~649~602 \\
12 & 0.119~759~673~589~890 & 0.115~129~254~649~358 & 0.115~129~254~649~772 \\
13 & 0.110~852~765~866~205 & 0.115~129~254~649~714 & 0.115~129~254~649~711 \\
14 & 0.119~099~529~780~039 & 0.115~129~254~649~727 & 0.115~129~254~649~700 \\
15 & 0.111~426~375~175~159 & 0.115~129~254~649~696 & 0.115~129~254~649~702 \\
16 & 0.118~596~725~346~641 & 0.115~129~254~649~702 & 0.115~129~254~649~702 \\
17 & 0.111~870~534~473~655 & 0.115~129~254~649~703 & 0.115~129~254~649~702 \\
18 & 0.118~201~666~639~655 & 0.115~129~254~649~702 & 0.115~129~254~649~702 \\
19 & 0.112~224~086~515~227 & 0.115~129~254~649~702 & 0.115~129~254~649~702 \\
20 & 0.117~883~505~752~715 & 0.115~129~254~649~702 & 0.115~129~254~649~702 \\
\tableline \rule{0pt}{3.5\jot}%
exact & 0.115~129~254~649~702 & 0.115~129~254~649~702 & 0.115~129~254~649~702
\end{tabular}
\end{table}

\begin{table}
\caption{Evaluation of $10^{-1}\,{\rm Li}_2 (0.99999)$ with the CNCT.}
\label{Tab_5_2}%
\begin{tabular}{lrrr}%
$n$%
& \multicolumn{1}{c}{${\bf S}_{n}$}%
& \multicolumn{1}{c}{$d_{n}^{(0)} \bigl(1, {\bf S}_0 \bigr)$}%
& \multicolumn{1}{c}{${\delta}_{n}^{(0)} \bigl(1, {\bf S}_0 \bigr)$}%
\\
\tableline \rule[-4pt]{0pt}{3.8\jot}%
 0 & 0.199~982~280~324~442 & 0.199~982~280~324~442 & 0.199~982~280~324~442 \\
 1 & 0.149~990~640~162~221 & 0.165~371~886~328~955 & 0.165~371~886~328~955 \\
 2 & 0.172~207~484~145~972 & 0.164~381~453~915~497 & 0.164~381~453~915~497 \\
 3 & 0.159~711~414~066~111 & 0.164~488~426~505~073 & 0.164~482~760~527~739 \\
 4 & 0.167~708~334~514~884 & 0.164~480~538~599~000 & 0.164~481~025~806~042 \\
 5 & 0.162~155~301~413~564 & 0.164~480~897~128~353 & 0.164~480~896~552~227 \\
 6 & 0.166~234~803~732~817 & 0.164~480~894~717~325 & 0.164~480~893~640~937 \\
 7 & 0.163~111~643~694~762 & 0.164~480~893~606~728 & 0.164~480~893~688~676 \\
 8 & 0.165~579~162~855~583 & 0.164~480~893~702~656 & 0.164~480~893~698~442 \\
 9 & 0.163~580~602~633~196 & 0.164~480~893~699~380 & 0.164~480~893~699~234 \\
10 & 0.165~232~198~806~363 & 0.164~480~893~699~272 & 0.164~480~893~699~288 \\
11 & 0.163~844~487~813~342 & 0.164~480~893~699~295 & 0.164~480~893~699~292 \\
12 & 0.165~026~841~360~043 & 0.164~480~893~699~292 & 0.164~480~893~699~293 \\
13 & 0.164~007~426~937~253 & 0.164~480~893~699~293 & 0.164~480~893~699~293 \\
14 & 0.164~895~394~970~447 & 0.164~480~893~699~293 & 0.164~480~893~699~293 \\
15 & 0.164~115~002~453~607 & 0.164~480~893~699~293 & 0.164~480~893~699~293 \\
\tableline \rule{0pt}{3.5\jot}%
exact & 0.164~480~893~699~293 & 0.164~480~893~699~293 & 0.164~480~893~699~293
\end{tabular}
\end{table}

\begin{table}
\caption{Evaluation of $10^{-1}\,{\rm Li}_3 (0.99999)$ with the CNCT.}
\label{Tab_5_3}%
\begin{tabular}{lrrr}%
$n$%
& \multicolumn{1}{c}{${\bf S}_{n}$}%
& \multicolumn{1}{c}{$d_{n}^{(0)} \bigl(1, {\bf S}_0 \bigr)$}%
& \multicolumn{1}{c}{${\delta}_{n}^{(0)} \bigl(1, {\bf S}_0\bigr)$}%
\\
\tableline \rule[-4pt]{0pt}{3.8\jot}%
 0 & 0.133~331~333~415~539 & 0.133~331~333~415~539 & 0.133~331~333~415~539 \\
 1 & 0.116~665~166~707~769 & 0.120~474~532~168~000 & 0.120~474~532~168~000 \\
 2 & 0.121~603~216~117~468 & 0.120~176~326~936~846 & 0.120~176~326~936~846 \\
 3 & 0.119~520~007~764~208 & 0.120~206~042~152~677 & 0.120~204~748~497~388 \\
 4 & 0.120~586~594~446~594 & 0.120~203~955~328~380 & 0.120~204~079~128~106 \\
 5 & 0.119~969~366~038~597 & 0.120~204~046~033~711 & 0.120~204~045~387~208 \\
 6 & 0.120~358~052~152~572 & 0.120~204~045~725~809 & 0.120~204~045~378~284 \\
 7 & 0.120~097~666~726~411 & 0.120~204~045~413~120 & 0.120~204~045~434~802 \\
 8 & 0.120~280~540~991~745 & 0.120~204~045~439~707 & 0.120~204~045~438~553 \\
 9 & 0.120~147~227~650~952 & 0.120~204~045~438~749 & 0.120~204~045~438~726 \\
10 & 0.120~247~386~435~540 & 0.120~204~045~438~729 & 0.120~204~045~438~733 \\
11 & 0.120~170~239~824~481 & 0.120~204~045~438~733 & 0.120~204~045~438~733 \\
12 & 0.120~230~916~813~857 & 0.120~204~045~438~733 & 0.120~204~045~438~733 \\
13 & 0.120~182~336~147~846 & 0.120~204~045~438~733 & 0.120~204~045~438~733 \\
14 & 0.120~221~833~436~597 & 0.120~204~045~438~733 & 0.120~204~045~438~733 \\
15 & 0.120~189~289~161~290 & 0.120~204~045~438~733 & 0.120~204~045~438~733 \\
\tableline \rule{0pt}{3.5\jot}%
exact & 0.120~204~045~438~733 & 0.120~204~045~438~733 & 0.120~204~045~438~733
\end{tabular}
\end{table}

\begin{table}
\caption{Evaluation of $10^4\,\Phi(0.99999,2,10000)$ with the CNCT.}%
\label{Tab_5_4}%
\begin{tabular}{lrrr}%
$n$%
& \multicolumn{1}{c}{${\bf S}_{n}$}%
& \multicolumn{1}{c}{$d_{n}^{(0)} \bigl(1, {\bf S}_0 \bigr)$}%
& \multicolumn{1}{c}{${\delta}_{n}^{(0)} \bigl(1, {\bf S}_0 \bigr)$}%
\\
\tableline \rule[-4pt]{0pt}{3.8\jot}%
 0 & 1.152~086~970~131~424 & 1.152~086~970~131~424 & 1.152~086~970~131~424 \\
 1 & 0.576~093~485~065~712 & 0.806~478~876~912~452 & 0.806~478~876~912~452 \\
 2 & 0.960~055~803~546~361 & 0.797~618~192~129~198 & 0.797~618~192~129~198 \\
 3 & 0.672~109~050~515~104 & 0.798~663~645~011~412 & 0.798~596~144~946~064 \\
 4 & 0.902~446~455~721~367 & 0.798~581~028~864~897 & 0.798~586~253~716~867 \\
 5 & 0.710~515~275~487~446 & 0.798~585~227~987~408 & 0.798~585~188~634~170 \\
 6 & 0.875~013~443~138~958 & 0.798~585~145~208~936 & 0.798~585~140~857~888 \\
 7 & 0.731~090~035~137~739 & 0.798~585~138~553~527 & 0.798~585~139~249~075 \\
 8 & 0.859~010~851~725~411 & 0.798~585~139~276~063 & 0.798~585~139~237~667 \\
 9 & 0.743~892~107~147~896 & 0.798~585~139~218~618 & 0.798~585~139~229~500 \\
10 & 0.848~536~431~310~159 & 0.798~585~139~219~617 & 0.798~585~139~222~908 \\
11 & 0.752~620~788~733~222 & 0.798~585~139~223~310 & 0.798~585~139~222~175 \\
12 & 0.841~150~625~351~432 & 0.798~585~139~222~720 & 0.798~585~139~222~491 \\
13 & 0.758~951~478~583~304 & 0.798~585~139~222~444 & 0.798~585~139~222~559 \\
14 & 0.835~664~028~102~224 & 0.798~585~139~222~555 & 0.798~585~139~222~550 \\
15 & 0.763~752~250~680~046 & 0.798~585~139~222~555 & 0.798~585~139~222~548 \\
16 & 0.831~428~053~244~539 & 0.798~585~139~222~546 & 0.798~585~139~222~548 \\
17 & 0.767~517~561~053~133 & 0.798~585~139~222~548 & 0.798~585~139~222~548 \\
18 & 0.828~059~092~035~324 & 0.798~585~139~222~548 & 0.798~585~139~222~548 \\
19 & 0.770~549~625~376~190 & 0.798~585~139~222~548 & 0.798~585~139~222~548 \\
20 & 0.825~315~796~529~872 & 0.798~585~139~222~548 & 0.798~585~139~222~548 \\
\tableline \rule{0pt}{3.5\jot}%
exact & 0.798~585~139~222~548 & 0.798~585~139~222~548 & 0.798~585~139~222~548
\end{tabular}
\end{table}

\begin{table}
\caption{Evaluation of $10^{-4}\,{_3} F_2(1,3/2,5;9/8,47/8;0.99999)$
with the CNCT.}%
\label{Tab_6_1}%
\begin{tabular}{lrrr}%
$n$%
& \multicolumn{1}{c}{${\bf S}_{n}$}%
& \multicolumn{1}{c}{$d_{n}^{(0)} \bigl(1, {\bf S}_0 \bigr)$}%
& \multicolumn{1}{c}{${\delta}_{n}^{(0)} \bigl(1, {\bf S}_0 \bigr)$}%
\\
\tableline \rule[-4pt]{0pt}{3.8\jot}%
 0 & 0.343~961~195~195~881 & 0.343~961~195~195~881 & 0.343~961~195~195~881 \\
 1 & 0.172~030~597~597~940 & 0.240~789~533~227~428 & 0.240~789~533~227~428 \\
 2 & 0.286~614~046~845~766 & 0.238~145~631~122~015 & 0.238~145~631~122~015 \\
 3 & 0.200~705~485~068~072 & 0.238~457~646~856~530 & 0.238~437~505~168~542 \\
 4 & 0.269~409~131~416~317 & 0.238~433~043~073~649 & 0.238~434~595~265~258 \\
 5 & 0.212~175~660~263~795 & 0.238~434~334~888~202 & 0.238~434~314~202~617 \\
 6 & 0.261~216~282~820~282 & 0.238~434~307~597~861 & 0.238~434~305~531~575 \\
 7 & 0.218~319~995~918~563 & 0.238~434~299~034~546 & 0.238~434~301~380~046 \\
 8 & 0.256~437~590~806~900 & 0.238~434~297~734~183 & 0.238~434~298~863~521 \\
 9 & 0.222~142~795~720~768 & 0.238~434~298~685~404 & 0.238~434~298~539~687 \\
10 & 0.253~309~971~367~287 & 0.238~434~298~885~696 & 0.238~434~298~711~196 \\
11 & 0.224~749~013~245~827 & 0.238~434~298~760~156 & 0.238~434~298~769~407 \\
12 & 0.251~104~853~004~796 & 0.238~434~298~751~345 & 0.238~434~298~766~552 \\
13 & 0.226~638~970~976~291 & 0.238~434~298~765~877 & 0.238~434~298~763~311 \\
14 & 0.249~467~024~022~419 & 0.238~434~298~763~853 & 0.238~434~298~763~230 \\
15 & 0.228~071~952~328~669 & 0.238~434~298~762~980 & 0.238~434~298~763~331 \\
16 & 0.248~202~731~012~486 & 0.238~434~298~763~382 & 0.238~434~298~763~332 \\
17 & 0.229~195~680~986~673 & 0.238~434~298~763~343 & 0.238~434~298~763~330 \\
18 & 0.247~197~365~060~190 & 0.238~434~298~763~322 & 0.238~434~298~763~330 \\
19 & 0.230~100~443~575~937 & 0.238~434~298~763~332 & 0.238~434~298~763~330 \\
20 & 0.246~378~830~994~286 & 0.238~434~298~763~330 & 0.238~434~298~763~330 \\
\tableline \rule{0pt}{3.5\jot}%
exact & 0.238~434~298~763~330 & 0.238~434~298~763~330 & 0.238~434~298~763~330
\end{tabular}
\end{table}

\begin{table}
\caption{Evaluation of $10^{-1}\,{_3} F_2(1,3,7;5/2,14;0.99999)$ with the CNCT.}%
\label{Tab_6_2}%
\begin{tabular}{lrrr}%
$n$%
& \multicolumn{1}{c}{${\bf S}_{n}$}%
& \multicolumn{1}{c}{$d_{n}^{(0)} \bigl(1, {\bf S}_0 \bigr)$}%
& \multicolumn{1}{c}{${\delta}_{n}^{(0)} \bigl(1, {\bf S}_0 \bigr)$}%
\\
\tableline \rule[-4pt]{0pt}{3.8\jot}%
 0 & 0.354~205~299~194~014 & 0.354~205~299~194~014 & 0.354~205~299~194~014 \\
 1 & 0.227~102~649~597~007 & 0.268~438~401~594~236 & 0.268~438~401~594~236 \\
 2 & 0.288~360~357~978~268 & 0.266~943~489~834~494 & 0.266~943~489~834~494 \\
 3 & 0.254~808~733~179~764 & 0.267~121~290~068~426 & 0.267~112~224~310~036 \\
 4 & 0.274~632~798~536~574 & 0.267~100~105~730~057 & 0.267~101~775~442~210 \\
 5 & 0.262~289~292~919~201 & 0.267~103~299~447~369 & 0.267~102~948~107~378 \\
 6 & 0.270~285~887~617~646 & 0.267~102~742~191~670 & 0.267~102~815~189~177 \\
 7 & 0.264~938~303~793~252 & 0.267~102~836~668~825 & 0.267~102~824~285~564 \\
 8 & 0.268~610~033~794~438 & 0.267~102~822~243~079 & 0.267~102~823~985~155 \\
 9 & 0.266~031~550~394~689 & 0.267~102~824~198~445 & 0.267~102~823~985~352 \\
10 & 0.267~878~112~865~182 & 0.267~102~823~960~840 & 0.267~102~823~984~751 \\
11 & 0.266~532~664~737~038 & 0.267~102~823~987~269 & 0.267~102~823~984~758 \\
12 & 0.267~528~211~474~989 & 0.267~102~823~984~510 & 0.267~102~823~984~761 \\
13 & 0.266~781~286~870~185 & 0.267~102~823~984~786 & 0.267~102~823~984~762 \\
14 & 0.267~348~762~560~130 & 0.267~102~823~984~759 & 0.267~102~823~984~762 \\
15 & 0.266~912~657~747~879 & 0.267~102~823~984~762 & 0.267~102~823~984~762 \\
16 & 0.267~251~339~082~039 & 0.267~102~823~984~762 & 0.267~102~823~984~761 \\
17 & 0.266~985~765~115~182 & 0.267~102~823~984~761 & 0.267~102~823~984~762 \\
18 & 0.267~195~878~973~722 & 0.267~102~823~984~762 & 0.267~102~823~984~762 \\
19 & 0.267~028~262~510~152 & 0.267~102~823~984~762 & 0.267~102~823~984~762 \\
20 & 0.267~163~009~854~165 & 0.267~102~823~984~762 & 0.267~102~823~984~762 \\
\tableline \rule{0pt}{3.5\jot}%
exact & 0.267~102~823~984~762 & 0.267~102~823~984~762 & 0.267~102~823~984~762
\end{tabular}
\end{table}

\begin{table}
\caption{Evaluation of $10^{-1}\,{_3} F_2(1,3,7;5/2,14;1)$ with the CNCT.}%
\label{Tab_6_3}%
\begin{tabular}{lrrr}%
$n$%
& \multicolumn{1}{c}{${\bf S}_{n}$}%
& \multicolumn{1}{c}{$d_{n}^{(0)} \bigl(1, {\bf S}_0 \bigr)$}%
& \multicolumn{1}{c}{${\delta}_{n}^{(0)} \bigl(1, {\bf S}_0 \bigr)$}%
\\
\tableline \rule[-4pt]{0pt}{3.8\jot}%
 0 & 0.354~212~896~979~703 & 0.354~212~896~979~703 & 0.354~212~896~979~703 \\
 1 & 0.227~106~448~489~852 & 0.268~443~680~394~043 & 0.268~443~680~394~043 \\
 2 & 0.288~366~524~620~961 & 0.266~948~705~538~902 & 0.266~948~705~538~902 \\
 3 & 0.254~813~300~376~035 & 0.267~126~514~679~686 & 0.267~117~448~402~341 \\
 4 & 0.274~638~493~612~452 & 0.267~105~329~111~381 & 0.267~106~998~932~606 \\
 5 & 0.262~294~169~832~612 & 0.267~108~523~033~192 & 0.267~108~171~668~587 \\
 6 & 0.270~291~370~710~880 & 0.267~107~965~739~446 & 0.267~108~038~742~291 \\
 7 & 0.264~943~330~016~988 & 0.267~108~060~223~478 & 0.267~108~047~839~255 \\
 8 & 0.268~615~409~392~370 & 0.267~108~045~796~597 & 0.267~108~047~538~821 \\
 9 & 0.266~036~655~402~122 & 0.267~108~047~752~131 & 0.267~108~047~539~018 \\
10 & 0.267~883~429~838~931 & 0.267~108~047~514~503 & 0.267~108~047~538~417 \\
11 & 0.266~537~813~952~027 & 0.267~108~047~540~935 & 0.267~108~047~538~424 \\
12 & 0.267~533~494~713~522 & 0.267~108~047~538~176 & 0.267~108~047~538~427 \\
13 & 0.266~786~462~107~999 & 0.267~108~047~538~452 & 0.267~108~047~538~428 \\
14 & 0.267~354~025~525~918 & 0.267~108~047~538~425 & 0.267~108~047~538~428 \\
15 & 0.266~917~848~923~282 & 0.267~108~047~538~428 & 0.267~108~047~538~428 \\
16 & 0.267~256~589~412~242 & 0.267~108~047~538~428 & 0.267~108~047~538~428 \\
17 & 0.266~990~966~387~224 & 0.267~108~047~538~428 & 0.267~108~047~538~428 \\
18 & 0.267~201~121~176~708 & 0.267~108~047~538~428 & 0.267~108~047~538~428 \\
19 & 0.267~033~470~369~207 & 0.267~108~047~538~428 & 0.267~108~047~538~428 \\
20 & 0.267~168~246~683~931 & 0.267~108~047~538~428 & 0.267~108~047~538~428 \\
\tableline \rule{0pt}{3.5\jot}%
exact & 0.267~108~047~538~428 & 0.267~108~047~538~428 & 0.267~108~047~538~428
\end{tabular}
\end{table}

\begin{table}
\caption{Evaluation of $10^{-5}\,\sum_{l=0}^{\infty} (2\,l+1)\,
  j_l(i\,0.9999\times0.7)\,h^{(1)}_l(i\,0.7)$ with the CNCT.}
\label{Tab_7_1}
\begin{tabular}{lrrr}%
$n$%
& \multicolumn{1}{c}{${\bf S}_{n}$}%
& \multicolumn{1}{c}{$d_{n}^{(0)} \bigl(1, {\bf S}_0\bigr)$}%
& \multicolumn{1}{c}{${\delta}_{n}^{(0)} \bigl(1, {\bf S}_0\bigr)$}%
\\
\tableline \rule[-4pt]{0pt}{3.8\jot}%
 0 & $-0.206~084~520~894~668$ & $-0.206~084~520~894~668$ & $-0.206~084~520~894~668$ \\
 1 & $-0.103~046~104~279~554$ & $-0.144~259~660~091~669$ & $-0.144~259~660~091~669$ \\
 2 & $-0.171~733~352~076~042$ & $-0.142~674~251~704~499$ & $-0.142~674~251~704~499$ \\
 3 & $-0.120~220~473~724~774$ & $-0.142~861~066~165~942$ & $-0.142~849~004~178~547$ \\
 4 & $-0.161~427~969~242~195$ & $-0.142~846~281~288~224$ & $-0.142~847~217~919~030$ \\
 5 & $-0.127~091~187~739~732$ & $-0.142~847~117~647~583$ & $-0.142~847~093~177~734$ \\
 6 & $-0.156~520~814~934~640$ & $-0.142~847~153~435~206$ & $-0.142~847~135~794~481$ \\
 7 & $-0.130~771~367~274~807$ & $-0.142~847~152~002~732$ & $-0.142~847~148~838~958$ \\
 8 & $-0.153~658~281~433~404$ & $-0.142~847~142~941~523$ & $-0.142~847~145~983~994$ \\
 9 & $-0.133~061~584~756~825$ & $-0.142~847~142~048~135$ & $-0.142~847~143~380~152$ \\
10 & $-0.151~784~408~324~881$ & $-0.142~847~143~286~397$ & $-0.142~847~143~026~940$ \\
11 & $-0.134~623~097~855~864$ & $-0.142~847~143~324~758$ & $-0.142~847~143~169~466$ \\
12 & $-0.150~463~209~081~729$ & $-0.142~847~143~181~923$ & $-0.142~847~143~211~999$ \\
13 & $-0.135~755~492~037~282$ & $-0.142~847~143~200~611$ & $-0.142~847~143~208~928$ \\
14 & $-0.149~481~826~233~342$ & $-0.142~847~143~210~852$ & $-0.142~847~143~207~036$ \\
15 & $-0.136~614~208~810~295$ & $-0.142~847~143~206~780$ & $-0.142~847~143~207~092$ \\
16 & $-0.148~724~109~578~581$ & $-0.142~847~143~206~921$ & $-0.142~847~143~207~139$ \\
17 & $-0.137~287~765~215~359$ & $-0.142~847~143~207~226$ & $-0.142~847~143~207~135$ \\
18 & $-0.148~121~427~644~456$ & $-0.142~847~143~207~123$ & $-0.142~847~143~207~135$ \\
19 & $-0.137~830~196~213~572$ & $-0.142~847~143~207~132$ & $-0.142~847~143~207~135$ \\
20 & $-0.147~630~649~333~917$ & $-0.142~847~143~207~137$ & $-0.142~847~143~207~135$ \\
21 & $-0.138~276~357~305~065$ & $-0.142~847~143~207~134$ & $-0.142~847~143~207~135$ \\
22 & $-0.147~223~291~777~104$ & $-0.142~847~143~207~135$ & $-0.142~847~143~207~135$ \\
23 & $-0.138~649~758~251~940$ & $-0.142~847~143~207~135$ & $-0.142~847~143~207~135$ \\
24 & $-0.146~879~773~958~865$ & $-0.142~847~143~207~135$ & $-0.142~847~143~207~135$ \\
25 & $-0.138~966~841~391~919$ & $-0.142~847~143~207~135$ & $-0.142~847~143~207~135$ \\
\tableline \rule{0pt}{3.5\jot}%
exact & $-0.142~847~143~207~135$ & $-0.142~847~143~207~135$ & $-0.142~847~143~207~135$
\end{tabular}
\end{table}

\end{document}